\documentclass[12pt]{article}
\usepackage{amsmath,amsthm,amssymb,dsfont,graphicx,xspace,epsfig,xcolor}
\usepackage[plain]{fullpage}
\usepackage{thm-restate}
\usepackage{hyperref}
\usepackage{authblk}
\usepackage{enumitem}
\usepackage{stmaryrd}
\usepackage{tikz, tkz-graph, tkz-berge}
\usetikzlibrary{calc}
\usetikzlibrary{decorations.pathmorphing}
\usetikzlibrary{decorations.pathreplacing}

\definecolor{g-green}{rgb}{0.235, 0.659, 0.322}
\definecolor{g-blue}{rgb}{0.0, 0.5, 1.0}
\usepackage{xcolor}

\newtheorem{theorem}{Theorem}
\newtheorem{corollary}[theorem]{Corollary}
\newtheorem{proposition}[theorem]{Proposition}
\newtheorem{lemma}[theorem]{Lemma}

\theoremstyle{definition}

\newtheorem{conjecture}[theorem]{Conjecture}

\newenvironment{proofclaim}[1][{\it Proof of claim. \hspace{0.066cm}}]
	{\noindent {}{#1}{}}{ \strut\hfill $\lozenge$\vspace{2ex}}

\newcommand{\dic}{\vec{\chi}}
\newcommand{\bid}{\overleftrightarrow}

\renewcommand{\leq}{\leqslant}
\renewcommand{\geq}{\geqslant}
\renewcommand{\phi}{\varphi}
\renewcommand{\epsilon}{\varepsilon}

\title{On the minimum number of arcs in $4$-dicritical oriented graphs \footnote{A 12-page extended abstract of this paper has been published in the proceedings of WG 2023 \cite{WG2023}. 

Research supported by research grant DIGRAPHS ANR-19-CE48-0013 and by the French government, through the EUR DS4H Investments in the Future project managed by the National Research Agency (ANR) with the reference number ANR-17-EURE-0004.
}}

\author[1]{Fr\'ed\'eric Havet}
\author[1]{Lucas Picasarri-Arrieta}
\author[1,2]{Cl\'ement Rambaud}

\affil[1]{Universit\'e C\^ote d'Azur, CNRS, Inria, I3S, Sophia Antipolis, France}
\affil[2]{DIENS, \'Ecole Normale Sup\'erieure, CNRS, PSL University, Paris, France}

\date{}

\begin{document}

\maketitle
\vspace{-10mm}
\begin{center}
{\small 
\texttt{$\{$frederic.havet, lucas.picasarri-arrieta, clement.rambaud$\}$@inria.fr}\\ 
}
\end{center}

\begin{abstract}
The dichromatic number $\dic(D)$ of a digraph $D$ is the minimum number of colours needed to colour the vertices of a digraph such that each colour class induces an acyclic subdigraph. 
A digraph $D$ is $k$-dicritical if $\dic(D) = k$ and each proper subdigraph $H$ of $D$ satisfies $\dic(H) < k$.

For integers $k$ and $n$, we define $d_k(n)$ (respectively $o_k(n)$) as the minimum number of arcs possible in a $k$-dicritical digraph (respectively oriented graph). 
Kostochka and Stiebitz have shown~\cite{kostochka_minimum_2020} that $d_4(n) \geq \frac{10}{3}n -\frac{4}{3}$. 
They also conjectured that there is a constant $c$ such that $o_k(n) \geq cd_k(n)$ for $k\geq 3$ and $n$ large enough.
This conjecture is known to be true for $k=3$ (Aboulker et al.~\cite{aboulker3dicritical}).

In this work, we prove that every $4$-dicritical oriented graph on $n$ vertices has at least $(\frac{10}{3}+\frac{1}{51})n-1$ arcs, showing the conjecture for $k=4$. 
We also characterise exactly the $4$-dicritical digraphs on $n$ vertices with exactly $\frac{10}{3}n -\frac{4}{3}$ arcs.
\end{abstract}

\section{Introduction}

Let $G$ be a graph. We denote by $V(G)$ its vertex set and by $E(G)$ its edge set; we set $n(G)=|V(G)|$ and $m(G)=|E(G)|$.
A {\bf $k$-colouring} of $G$ is a function $\phi:V(G) \to [k]$. It is {\bf proper} if for every edge $uv \in E(G)$, $\phi(u) \neq \phi(v)$.
The smallest integer $k$ such that $G$ has a proper $k$-colouring is the {\bf chromatic number}, and is denoted by $\chi(G)$.
Since $\chi$ is non-decreasing with respect to the subgraph relation, it is natural to consider the minimal graphs (for this relation) which are not $(k-1)$-colourable.
Following this idea, Dirac defined {\bf $k$-critical} graphs as the graphs $G$ with $\chi(G) = k$ and $\chi(H) < k$ for every proper subgraph $H$ of $G$.
A first property of $k$-critical graph is that their minimum degree is at least $k-1$. Indeed, if a vertex $v$ has degree at most $k-2$, then a $(k-1)$-colouring of $G-v$ can be easily extended to $G$, contradicting the fact that $\chi(G)=k$.
As a consequence, the number of edges in a $k$-critical graph is at least $\frac{k-1}{2}n$.
This bound is tight for complete graphs and odd cycles, but Dirac~\cite{dirac1957theorem} proved an inequality of the form $m \geq \frac{k-1+\epsilon_k}{2}n - c_k$ for every $n$-vertex $k$-critical graph with $m$ edges, for some $c_k$ and $\epsilon_k>0$. This shows that, for $n$ sufficiently large, the average degree of a $k$-critical graph is at least $k-1+\epsilon_k$.
This initiated the quest after the best lower bound on the number of edges in $n$-vertex $k$-critical graphs. This problem was almost completely solved by Kostochka and Yancey in 2014~\cite{kostochka_ores_2014}.

\begin{theorem}[Kostochka and Yancey~\cite{kostochka_ores_2014}]\label{theorem:kostochka_and_yancey_general_result} ~\\
Every $k$-critical graph on $n$ vertices has at least $\frac{1}{2}(k-\frac{2}{k-1})n - \frac{k(k-3)}{2(k-1)}$ edges.
For every $k$, this bound is tight for infinitely many values of $n$.
\end{theorem}
Kostochka and Yancey~\cite{kostochka2018brooks} also characterised $k$-critical graphs for which this inequality is an equality,
and all of them contain a copy of $K_{k-2}$, the complete graph on $k-2$ vertices.
This motivated the following conjecture of Postle~\cite{Postle-EJC17}.

\begin{conjecture}[Postle~\cite{Postle-EJC17}]\label{conjecture:postle}
    For every integer $k \geq 4$, there exists $\epsilon_k > 0$ such that every $k$-critical $K_{k-2}$-free graph $G$
    on $n$ vertices has at least $\frac{1}{2}\left(k-\frac{2}{k-1} + \epsilon_k\right)n - \frac{k(k-3)}{2(k-1)}$ edges.
\end{conjecture}

For $k=4$, the conjecture trivially holds as there is no $K_2$-free $4$-critical graph.
Moreover, this conjecture has been confirmed for $k=5$ by Postle~\cite{Postle-EJC17}, for $k=6$ by Gao and Postle~\cite{gao_minimal_2018}, and for $k\geq 33$ by Gould, Larsen, and Postle~\cite{gould_structure_2021}.

\medskip

Let $D$ be a digraph. We denote by $V(D)$ its vertex set and by $A(D)$ its arc set;  we set $n(D)=|V(D)|$ and $m(D)=|A(D)|$. 
A {\bf $k$-colouring} of $D$ is a function $\phi:V(D) \to [k]$. It is a {\bf $k$-dicolouring} if every directed cycle $C$ in $D$ is not monochromatic for $\phi$ (that is $|\phi(V(C))| > 1$). Equivalently, it is a $k$-dicolouring  if every colour class induces an acyclic subdigraph.
The smallest integer $k$ such that $D$ has a $k$-dicolouring is the {\bf dichromatic number} of $D$ and is denoted by $\dic(D)$.

A {\bf digon} in $D$ is a pair of opposite arcs between two vertices. Such a pair of arcs $\{uv,vu\}$ is denoted by $[u,v]$. 
We say that $D$ is a {\bf bidirected graph} if every pair of adjacent vertices forms a digon.
In this case, $D$ can be viewed as obtained from an undirected graph $G$ by replacing each edge $\{u,v\}$ of $G$ by the digon $[u,v]$. 
We say that $D$ is a bidirected $G$, and we denote it by $\bid{G}$.
Observe that $\chi(G) = \dic(\bid{G})$. Thus every statement on proper colouring of undirected graphs can be seen as a statement on dicolouring of bidirected graphs.

Exactly as in the undirected case, one can define {\bf $k$-dicritical} digraphs to be digraphs $D$ with $\dic(D)=k$ and
$\dic(H)<k$ for every proper subdigraph $H$ of $D$.
It is easy to check that if $G$ is a $k$-critical graph, then $\bid{G}$ is $k$-dicritical.
Kostochka and Stiebitz~\cite{kostochka_minimum_2020} conjectured  that the $k$-dicritical digraphs with the minimum number of arcs are bidirected graphs. Thus they conjectured the following generalisation of Theorem~\ref{theorem:kostochka_and_yancey_general_result} to digraphs.

\begin{conjecture}[Kostochka and Stiebitz~\cite{kostochka_minimum_2020}]\label{conj:KosSti}
    Let $k \geq 2$.
    Every $k$-dicritical digraph on $n$ vertices has at least $(k-\frac{2}{k-1})n - \frac{k(k-3)}{k-1}$ arcs.
    Moreover, equality holds only if $D$ is bidirected.
\end{conjecture}
In the case $k=2$, this conjecture is easy and weak as it states that a $2$-dicritical digraph
on $n$ vertices has at least two arcs, while, for all $n\geq 2$, the unique $2$-dicritical digraph of order $n$ is the directed $n$-cycle which has $n$ arcs.
The case $k=3$ of the conjecture has been confirmed by Kostochka and Stiebitz~\cite{kostochka_minimum_2020}.
Using a Brooks-type result for digraphs due to 
Harutyunyan and Mohar~\cite{harutyunyan2011gallai}, they proved the following: if $D$ is a $3$-dicritical digraph of order $n\geq 3$, then $m(D)\geq 2n$ and equality holds if and only if $n$ is odd and $D$ is a bidirected odd cycle.
The conjecture has also been proved for $k=4$ by Kostochka and Stiebitz~\cite{kostochka_minimum_2020}.
However, the conjecture is open for every $k\geq 5$.
Recently, this problem has been investigated by Aboulker and Vermande~\cite{aboulker2022various} who proved the weaker bound
$(k-\frac{1}{2}-\frac{2}{k-1})n - \frac{k(k-3)}{k-1}$ for the number of arcs in an $n$-vertex $k$-dicritical digraph.

For integers $k$ and $n$, let $d_k(n)$ denote the minimum number of arcs in a $k$-dicritical digraph of order $n$. 
By the above observations, $d_2(n) = n$ for all $n\geq 2$, and $d_3(n) \geq 2n$ for all possible $n$, and equality holds if and only if $n$ is odd and $n \geq 3$. Moreover, if $n$ is even then $d_3(n)=2n+1$ (see~\cite{aboulker3dicritical}). 
\medskip

Kostochka and Stiebitz~\cite{kostochka_number_2000} showed that if a $k$-critical graph $G$
 is triangle-free (that is has no cycle of length $3$), then $m(G)/n(G) \geq k - o(k)$
 as $k \to + \infty$. Informally, this means that
the minimum average degree of a $k$-critical triangle-free graph is (asymptotically) twice the minimum average degree of a $k$-critical graph.
Similarly to this undirected case, it is expected that the  minimum number of arcs in a $k$-dicritical digraph
of order $n$ is larger than $d_k(n)$ if we impose this digraph to have no short directed cycles, and in particular if the digraph is an {\bf oriented graph}, that is a digraph with no digon.
Let $o_k(n)$ denote the minimum number of arcs in a $k$-dicritical oriented graph of order $n$ (with the convention $o_k(n)=+\infty$ if there is no $k$-dicritical
oriented graph of order $n$). 
Clearly $o_k(n) \geq d_k(n)$.

\begin{conjecture}[Kostochka and Stiebitz~\cite{kostochka_minimum_2020}]\label{conj:alpha}
For any $k\geq 3$, there is a constant $\alpha_k>0$ such that $o_k(n) > (1+\alpha_k) d_k(n)$ for $n$ sufficiently large.
\end{conjecture}

For $k=3$, this conjecture has been recently confirmed by Aboulker, Bellitto, Havet, and Rambaud~\cite{aboulker3dicritical}
who proved that $o_3(n) \geq (2+\frac{1}{3})n + \frac{2}{3}$.

In view of Conjecture~\ref{conjecture:postle}, Conjecture~\ref{conj:alpha} can be generalized to $\bid{K_{k-2}}$-free digraphs.

\begin{conjecture}\label{conj:beta}
For any $k\geq 4$, there is a constant $\beta_k>0$ such that every $k$-dicritical $\bid{K_{k-2}}$-free digraph $D$  on $n$ vertices has at least $(1 + \beta_k)d_k(n)$ arcs.
\end{conjecture}

Together with Conjecture~\ref{conj:KosSti}, this conjecture would imply the following generalisation of Conjecture~\ref{conjecture:postle}.

\begin{conjecture}\label{conj:dicritical}
    For every integer $k \geq 4$, there exists $\epsilon_k > 0$ such that every $k$-dicritical $\bid{K_{k-2}}$-free digraph $D$
    on $n$ vertices has at least $(k-\frac{2}{k-1} + \epsilon_k)n - \frac{k(k-3)}{k-1}$ arcs.
\end{conjecture}

A $\bid{K_{2}}$-free digraph is an oriented graph, and there are infinitely many $4$-dicritical oriented graphs. Thus, while Conjecture~\ref{conjecture:postle} holds vacuously for $k=4$, this is not the case for Conjecture~\ref{conj:dicritical}.
In this paper, we prove that Conjectures~\ref{conj:alpha}, \ref{conj:beta}, and~\ref{conj:dicritical} hold for $k=4$.

\begin{theorem}\label{thm:main_corollary_oriented_4_dicritical}
    If $\vec{G}$ is a $4$-dicritical oriented graph, then 
    \[\displaystyle m(\vec{G}) \geq \left(\frac{10}{3}+\frac{1}{51}\right)n(\vec{G})-1. \]
\end{theorem}

To prove Theorem~\ref{thm:main_corollary_oriented_4_dicritical}, we use an approach similar to the proof of the case $k=5$ of Conjecture~\ref{conjecture:postle} by Postle~\cite{Postle-EJC17}. This proof is based on the potential method, which was first popularised by Kostochka and Yancey~\cite{kostochka_ores_2014} when they proved Theorem~\ref{theorem:kostochka_and_yancey_general_result}.
The idea is to prove a more general result on every $4$-dicritical digraphs that takes into account the digons. 

With a slight abuse, we call {\bf digon} a subdigraph isomorphic to $\bid{K_2}$, the bidirected complete graph on two vertices. We also call {\bf bidirected triangle} a subdigraph isomorphic to $\bid{K_3}$, the bidirected complete graph on three vertices. 
A {\bf packing} of digons and bidirected triangles is a set of vertex-disjoint digons and bidirected triangles.
To take into account the digons, we define a parameter $T(D)$
as follows.
\begin{equation*}
T(D)  =  \max\{ d+2t \mid \mbox{there exists a packing of $d$ digons and $t$ bidirected triangles}\}
\end{equation*}
Clearly, $T(D)=0$ if and only if $D$ is an oriented graph.

Let $\epsilon,\delta$ be fixed non-negative real numbers.
We define the {\bf potential} (with respect to $\epsilon$ and $\delta$) of a digraph $D$ to be \[\rho(D) = \left(\frac{10}{3}+\epsilon\right)n(D)-m(D)-\delta T(D).\]

Thus Theorem~\ref{thm:main_corollary_oriented_4_dicritical} can be rephrased as follows. 
\addtocounter{theorem}{-1}
\begin{theorem}
Set $\epsilon=\frac{1}{51}$ and $\delta = 6\epsilon = \frac{2}{17}$.
If $\vec{G}$ is a $4$-dicritical oriented graph, then 
$\rho(\vec{G}) \leq 1$. 
\end{theorem}

In fact, we prove a more general statement which holds for every $4$-dicritical digraph (with or without digons), except for some exceptions called the {\bf $4$-Ore digraphs}.
Those digraphs, which are formally defined in Section~\ref{sec:prop_4ore}, are the bidirected graphs whose underlying graph is one of the $4$-critical graphs reaching equality in Theorem~\ref{theorem:kostochka_and_yancey_general_result}. In particular, every $4$-Ore digraph $D$ has 
$\frac{10}{3}n(D) - \frac{4}{3}$ arcs.
Moreover, the statement holds for all non-negative constants $\epsilon$ and $\delta$
satisfying the following inequalities:
\begin{itemize}
    \item $\delta \geq 6\epsilon$;
    \item $3\delta - \epsilon \leq \frac{1}{3}$;
\end{itemize}

\begin{theorem}\label{thm:main_thm_potential}
Let $\epsilon,\delta \geq 0$ be constants satisfying the aforementioned inequalities.
If $D$ is a $4$-dicritical digraph with $n$ vertices, then
\begin{enumerate}[label=(\roman*)]
    \item $\rho(D) \leq \frac{4}{3} + \epsilon n - \delta \frac{2(n-1)}{3}$
        if $D$ is $4$-Ore, and
    \item $\rho(D) \leq 1$ otherwise.
\end{enumerate}
\end{theorem}

In order to provide some intuition to the reader, let us briefly describe the main ideas of our proof. We will consider a minimum counterexample $D$ to Theorem~\ref{thm:main_thm_potential}, and  show that every subdigraph of $D$ must have large potential. To do so, we need to construct some smaller 4-dicritical digraphs to leverage the minimality of $D$. These smaller 4-dicritical digraphs will be constructed by identifying some vertices of $D$. This is why, in the definition of the potential, we consider $T(D)$ instead of the number of digons: when identifying a set of vertices, the number of digons may be arbitrary larger in the resulting digraph, but $T(D)$ increases at most by 1. Using the fact that every subdigraph of $D$ has large potential, we will prove that some subdigraphs are forbidden in $D$. Using this, we get the final contradiction by a discharging argument.

\medskip

In addition to Theorem~\ref{thm:main_corollary_oriented_4_dicritical}, Theorem~\ref{thm:main_thm_potential} has also the following consequence when we take $\epsilon=\delta=0$.

\begin{corollary}
    If $D$ is a $4$-dicritical digraph, then $m(D) \geq \frac{10}{3}n(D) - \frac{4}{3}$.
    Moreover, equality holds if and only if $D$ is $4$-Ore, otherwise $m(D) \geq \frac{10}{3}n(D) - 1$.
\end{corollary}

This is a slight improvement on a result of Kostochka and Stiebitz~\cite{kostochka_minimum_2020} who proved the inequality $m(D) \geq \frac{10}{3}n(D)-\frac{4}{3}$ without characterising the equality case.

Another interesting consequence of our result is the following bound on the number of vertices in a $4$-dicritical oriented graph embedded on a fixed surface.
Since a graph on $n$ vertices embedded on a surface of Euler characteristic $c$ has at most $3n-3c$ edges, we immediately deduce the following from Theorem~\ref{thm:main_corollary_oriented_4_dicritical}.

\begin{corollary}
If $\vec{G}$ is a $4$-dicritical oriented graph embedded on a surface of
Euler characteristic $c$, then 
$n(\vec{G}) \leq \frac{17}{6}(1-3c)$.
\end{corollary}

The previous best upper bound was $n(\vec{G}) \leq  4 - 9c$~\cite{kostochka_minimum_2020}.

\medskip

In Section~\ref{sec:prop_4ore} we prove some first preliminary results on $4$-Ore digraphs, before proving Theorem~\ref{thm:main_thm_potential} in Section~\ref{sec:main_proof}. 
In Section~\ref{section:upperbounds}, we show that $o_k(n)\leq (2k-\frac{7}{2})n$ for every fixed $k$ and infinitely many values of $n$. The proof is strongly based on the proof of~\cite[Theorem 4.4]{aboulker3dicritical}, which shows $o_k(n)\leq (2k-3)n$ for every fixed $k,n$ (with $n$ large enough). For $k=4$, the construction implies in particular that there is a $4$-dicritical oriented graph with $76$ vertices and $330$ arcs, and there are infinitely many $4$-dicritical oriented graphs with $m/n \leq 9/2$.

\section{The \texorpdfstring{$4$}{4}-Ore digraphs and their properties}\label{sec:prop_4ore}

We start with a few notations.
We denote by $\llbracket x_1, \dots, x_n\rrbracket$ the bidirected path with vertex set $\{x_1, \dots, x_n\}$ in this order.
If $x_1=x_n$, $\llbracket x_1, \dots, x_n \rrbracket$ denotes the bidirected cycle of order $n$ with cyclic order $x_1, \dots, x_n$.
If $D$ is a digraph, for any $X \subseteq V(D)$, $D-X$ is the subdigraph induced by $V(D)\setminus X$. We abbreviate $D-\{x\}$ into $D-x$.
Moreover, for any $F\subseteq  V(D) \times V(D)$, $D\setminus F$ is the subdigraph $(V(D), A(D)\setminus F)$ and $D\cup F$ is the digraph $(V(D), A(D)\cup F)$

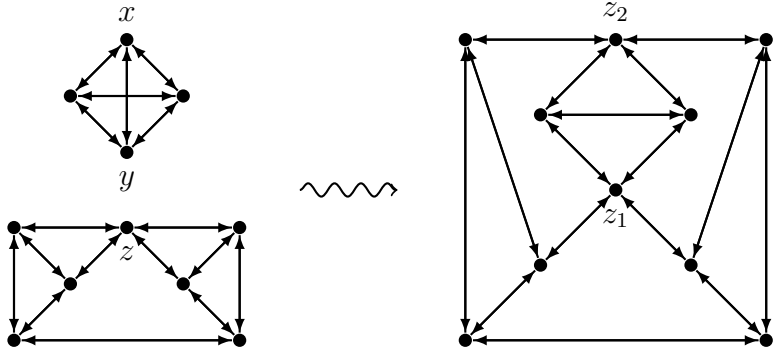
\begin{figure}
  \begin{minipage}{\linewidth}
    \begin{center}	
      \begin{tikzpicture}[thick,scale=1, every node/.style={transform shape}]
        \tikzset{vertex/.style = {circle,fill=black,minimum size=5pt,
                                        inner sep=0pt}}
        \tikzset{edge/.style = {->,> = latex}}
	  
        \node[vertex] (u1) at  (-7.25,-0.75) {};
        \node[vertex, label=above:$x$] (u2) at  (-6.5,0) {};
        \node[vertex] (u3) at  (-5.75,-0.75) {};
        \node[vertex, label=below:$y$] (u4) at  (-6.5,-1.5) {};
        
        \node[vertex, label=below:$z$] (u11) at  (-6.5,-2.5) {};
        
        \node[vertex] (u5) at  (-5,-2.5) {};
        \node[vertex] (u6) at  (-8,-2.5) {};
        
        \node[vertex] (u7) at  (-7.25,-3.25) {};
        \node[vertex] (u8) at  (-8,-4) {};
        \node[vertex] (u9) at  (-5.75,-3.25) {};
        \node[vertex] (u10) at  (-5,-4) {};
        \draw[edge] (u1) to (u2);
        \draw[edge] (u2) to (u1);
        \draw[edge] (u4) to (u2);
        \draw[edge] (u2) to (u4);
        \draw[edge] (u1) to (u3);
        \draw[edge] (u3) to (u1);
        \draw[edge] (u1) to (u4);
        \draw[edge] (u4) to (u1);
        \draw[edge] (u2) to (u3);
        \draw[edge] (u3) to (u2);
        \draw[edge] (u3) to (u4);
        \draw[edge] (u4) to (u3);
        \draw[edge] (u11) to (u5);
        \draw[edge] (u5) to (u11);
        \draw[edge] (u11) to (u6);
        \draw[edge] (u6) to (u11);
        \draw[edge] (u11) to (u7);
        \draw[edge] (u7) to (u11);
        \draw[edge] (u8) to (u7);
        \draw[edge] (u7) to (u8);
        \draw[edge] (u8) to (u10);
        \draw[edge] (u10) to (u8);
        \draw[edge] (u9) to (u10);
        \draw[edge] (u10) to (u9);
        \draw[edge] (u11) to (u9);
        \draw[edge] (u9) to (u11);
        \draw[edge] (u6) to (u7);
        \draw[edge] (u7) to (u6);
        \draw[edge] (u8) to (u6);
        \draw[edge] (u6) to (u8);
        \draw[edge] (u5) to (u9);
        \draw[edge] (u9) to (u5);
        \draw[edge] (u10) to (u5);
        \draw[edge] (u5) to (u10);
        \draw [->,decorate,decoration=snake] (-4.2,-2) -- (-2.9,-2);
        \node[vertex] (v1) at  (-1,-1) {};
        \node[vertex,label=above:$z_2$] (v2) at  (0,0) {};
        \node[vertex] (v3) at  (1,-1) {};
        \node[vertex,label=below:$z_1$] (v4) at  (0,-2) {};
        \node[vertex] (v5) at  (2,0) {};
        \node[vertex] (v6) at  (-2,0) {};
        \node[vertex] (v7) at  (-1,-3) {};
        \node[vertex] (v8) at  (-2,-4) {};
        \node[vertex] (v9) at  (1,-3) {};
        \node[vertex] (v10) at  (2,-4) {};
        \draw[edge] (v1) to (v2);
        \draw[edge] (v2) to (v1);
        \draw[edge] (v1) to (v3);
        \draw[edge] (v3) to (v1);
        \draw[edge] (v1) to (v4);
        \draw[edge] (v4) to (v1);
        \draw[edge] (v2) to (v3);
        \draw[edge] (v3) to (v2);
        \draw[edge] (v3) to (v4);
        \draw[edge] (v4) to (v3);
        \draw[edge] (v2) to (v5);
        \draw[edge] (v5) to (v2);
        \draw[edge] (v2) to (v6);
        \draw[edge] (v6) to (v2);
        \draw[edge] (v4) to (v7);
        \draw[edge] (v7) to (v4);
        \draw[edge] (v8) to (v7);
        \draw[edge] (v7) to (v8);
        \draw[edge] (v8) to (v10);
        \draw[edge] (v10) to (v8);
        \draw[edge] (v9) to (v10);
        \draw[edge] (v10) to (v9);
        \draw[edge] (v4) to (v9);
        \draw[edge] (v9) to (v4);
        \draw[edge] (v6) to (v7);
        \draw[edge] (v7) to (v6);
        \draw[edge] (v8) to (v6);
        \draw[edge] (v6) to (v8);
        \draw[edge] (v5) to (v9);
        \draw[edge] (v9) to (v5);
        \draw[edge] (v10) to (v5);
        \draw[edge] (v5) to (v10);
      \end{tikzpicture}
      \caption{An example of a 4-Ore digraph obtained by an Ore-composition of two smaller 4-Ore digraphs, with replaced digon $[x,y]$ and split vertex $z$.}
      \label{fig:4ore}
    \end{center}    
  \end{minipage}
\end{figure}

Let $D_1,D_2$ be two bidirected graphs, $[x,y] \subseteq A(D_1)$, and $z \in V(D_2)$.
An {\bf Ore-composition} $D$ of $D_1$ and $D_2$ with {\bf replaced digon $[x,y]$} and {\bf split vertex $z$} is a digraph obtained by removing $[x,y]$ of $D_1$ and $z$ of $D_2$,
and adding the set of arcs $\{xz_1 \mid zz_1\in A(D_2) \mbox{~and~} z_1\in Z_1\}$, $\{z_1x \mid z_1z\in A(D_2) \mbox{~and~} z_1\in Z_1\}$, $\{yz_2 \mid zz_2\in A(D_2) \mbox{~and~} z_2\in Z_2\}$, $\{z_2y \mid z_2z\in A(D_2) \mbox{~and~} z_2\in Z_2\}$, where $(Z_1,Z_2)$ is a partition of $N_{D_2}(z)$
into non-empty sets.
We call $D_1$ the {\bf digon side} and $D_2$ the {\bf split side} of the Ore-composition.
The class of the {\bf $4$-Ore digraphs} is the smallest class containing $\bid{K_4}$ which is stable under Ore-composition.
See Figure~\ref{fig:4ore} for an example of a 4-Ore digraph.
Observe that all the $4$-Ore-digraphs are bidirected. 

\begin{proposition}[Dirac~\cite{dirac1964structure}, see also~\cite{kostochka2018brooks}]\label{lemma:4ore_4_dicritical}
    $4$-Ore digraphs are $4$-dicritical.
\end{proposition}

\begin{proof}
    One can easily show that a bidirected digraph is $4$-dicritical if and only if its undirected underlying graph is $4$-critical.
    Then the result follows from the undirected analogous proved by~\cite{dirac1964structure}.
\end{proof}

\begin{restatable}{lemma}{orecompositiononlyinterestingifbidirected}\label{lemma:ore_composition_only_interesting_if_bidirected}

Let $D$ be a $4$-dicritical bidirected digraph and $v \in V(D)$.
    Let $(N_1^+,N_2^+)$  and $(N_1^-, N_2^-)$ be two partitions of $N(v)$.
    Consider $D'$ the digraph with vertex set $V(D) \setminus \{v\} \cup \{v_1,v_2\}$
    with $N^+(v_i)=N^+_i, N^-(v_i)=N^-_i$ for $i=1,2$ and $D'\langle V(D) \setminus \{v\} \rangle = D-v$.
    Then $D'$ has a $3$-dicolouring with $v_1$ and $v_2$ coloured the same except if $N^+_1=N^-_1$ (that is $D'$ is bidirected).
\end{restatable}
\begin{proof}
    Suppose that $D'$ is not bidirected. Consider a vertex $u \in N_D(v)$ such that $v_1u,uv_2 \in A(D')$ or $v_2u,uv_1 \in A(D')$.
    Without loss of generality, suppose $v_1u,uv_2 \in A(D')$.
    As $D$ is $4$-dicritical, $D \setminus [u,v]$ has a proper $3$-dicolouring $\phi$. We set $\phi(v_1)=\phi(v_2)=\phi(v)$ and claim that $\phi$ is a $3$-dicolouring of $D'$.
    To show that, observe that $\phi$ is a proper $3$-colouring of the underlying undirected graph of $D' \setminus \{v_1u,uv_2\}$, and so $\phi$ is a $3$-dicolouring of $D'$ as wanted.
\end{proof}

\begin{restatable}{lemma}{computationTminusonevertex}
\label{lemma:computation_T_minus_one_vertex}
Let $D$ be a digraph.
If $v$ is a vertex of $D$, then $T(D-v) \geq T(D)-1$.    
\end{restatable}
\begin{proof}
Let $M$ be a packing of $d$ digons and $t$ bidirected triangles in $H$ such that $d+2t=T(D)$.
If $v$ belongs to a digon $[u,v]$ in $M$, then $M \setminus \{[u,v]\}$ witnesses the fact that $T(D-v) \geq T(D)-1$.
If $v$ belongs to a bidirected triangle $\llbracket u,v,w,u \rrbracket$, then $M \setminus \{\llbracket u,v,w,u \rrbracket\} \cup [u,w]$ witnesses the fact that $T(D-v) \geq T(D) -2 + 1$.
Otherwise $T(D-v) \geq T(D)$.
\end{proof}

\begin{restatable}{lemma}{packingandorecomposition}
\label{lemma:packing_and_ore_composition}
    If $D_1,D_2$ are two digraphs, and $D$ is an Ore-composition of $D_1$ and
$D_2$, then $T(D) \geq T(D_1)+T(D_2)-2$. Moreover, if $D_1$ or $D_2$ is isomorphic to $\bid{K_4}$, then $T(D) \geq T(D_1)+T(D_2)-1$.
\end{restatable}
\begin{proof}
Let $D$ be the Ore-composition of $D_1$ (the digon side with replaced digon $[x,y]$) and $D_2$ (the split side with split vertex $z$).
One can easily see that $T(D) \geq T(D_1 - x) + T(D - z) \geq T(D_1)+T(D_2)-2$
by Lemma~\ref{lemma:computation_T_minus_one_vertex}.
Moreover, if $D_1$ (resp. $D_2$) is a copy of $\bid{K_4}$, then
$T(D_1-x)=2 = T(D_1)$ (resp. $T(D_2-z) = 2 = T(D_2)$) and therefore
$T(D) \geq T(D_1) + T(D_2)-1$.
\end{proof}

\begin{restatable}{lemma}{matchingore}
\label{lemma:matching_ore}
    If $D$ is $4$-Ore, then $T(D) \geq \frac{2}{3}(n(D)-1)$.
\end{restatable}
\begin{proof}
If $D$ is $\bid{K_4}$, then the result is clear.
Suppose now that $D$ is an Ore-composition of $D_1$ and $D_2$.
Then $n(D) = n(D_1) + n(D_2) - 1$ and, by
Lemma~\ref{lemma:packing_and_ore_composition},
$T(D) \geq T(D_1)+T(D_2)-2$.
By induction, $T(D_1) \geq \frac{2}{3}(n(D_1)-1)$ and $T(D_2) \geq \frac{2}{3}(n(D_2)-1)$,
and so $T(D) \geq \frac{2}{3}(n(D_1)+n(D_2)-1 -1) = \frac{2}{3}(n(D)-1)$.
\end{proof}

Let $D$ be a digraph.
A {\bf diamond} in $D$ is a subdigraph isomorphic to
$\bid{K_4}$ minus a digon $[u,v]$, with vertices different from $u$ and $v$
having degree $6$ in $D$.
An {\bf emerald} in $D$ is a subdigraph isomorphic to $\bid{K_3}$
whose vertices have degree $6$ in $D$.

Let $R$ be an induced subdigraph of $D$ with $n(R) < n(D)$.
The {\bf boundary} of $R$ in $D$, denoted by $\partial_D(R)$, or simply $\partial(R)$ when $D$ is clear from the context, 
is the set of vertices of $R$ having a neighbour in $V(D)\setminus R$. 
We say that $R$ is {\bf Ore-collapsible} if the boundary of $R$ contains exactly two vertices $u$ and $v$ and $R \cup [u,v]$ is $4$-Ore.

\begin{restatable}{lemma}{twofour}
\label{lemma:2.4}
    If $D$ is $4$-Ore and $v \in V(D)$, then there exists either an Ore-collapsible
subdigraph of $D$ disjoint from $v$ or an emerald of $D$ disjoint from $v$.
\end{restatable}
\begin{proof}
If $D$ is a copy of $\bid{K_4}$, then $D-v$ is an emerald.
Otherwise, $D$ is the Ore-composition of two $4$-Ore digraphs: $D_1$ the digon side with replaced digon $[x,y]$, and $D_2$ the split side with split vertex $z$.
If $v \in V(D_2-z)$, then $D_1$ is an Ore-collapsible subdigraph with boundary $\{x,y\}$. 
Otherwise $v \in V(D_1)$ and we apply induction on $D_2$ to find an emerald or an Ore-collapsible subdigraph in $D_2$ disjoint from $z$.
\end{proof}

\begin{restatable}{lemma}{twofive}
\label{lemma:2.5}
    If $D \neq \bid{K_4}$ is $4$-Ore and $T$ is a copy of $\bid{K_3}$ in $D$, then there exists either an Ore-collapsible subdigraph of $D$ disjoint from $T$ or
an emerald of $D$ disjoint from $T$.
\end{restatable}
\begin{proof}
As $D$ is not $\bid{K_4}$, it is an Ore-composition of two $4$-Ore digraphs: $D_1$ the digon side with replaced digon $[x,y]$, and $D_2$ the split side with split vertex $z$.
As $x$ and $y$ are not adjacent, we have either $T \subseteq D_1$, $T \subseteq D_2-z$, or $T$ contains a vertex $w \in \{x,y\}$ and two vertices in $V(D_2-z)$.

If $T \subseteq D_1$, then by Lemma~\ref{lemma:2.4}, in $D_2$ there exists either an Ore-collapsible subdigraph $O$ or an emerald $E$ disjoint from $z$. In the former case $O$ is an Ore-collapsible subdigraph of $D$ disjoint from $T$, and in the later one $E$ is an emerald in $D$ disjoint from $T$. 

If $T \subseteq D_2-z$, then $D_1\setminus \{x,y\}$ is an Ore-collapsible subdigraph disjoint from $T$.

Assume now that  $T$ contains a vertex $w \in \{x,y\}$ and two vertices in $V(D_2-z)$.
Without loss of generality, we may assume that $y \not\in T$.
Let $z_1$ and $z_2$ be the two vertices of $T$ disjoint from $w$. Then $\{z, z_1, z_2\}$ induces a bidirected triangle $T'$ in $D_2$.
If $D_2\neq \bid{K_4}$, then by induction in $D_2$, there exists either an Ore-collapsible subdigraph $O$ or an emerald $E$ disjoint from $T'$. In the former case $O$ is an Ore-collapsible subdigraph of $D$ disjoint from $T$, and in the later one $E$ is an emerald in $D$ disjoint from $T$.

Henceforth we may assume that $D_2=\bid{K_4}$.
This implies that $y$ has exactly one neighbour in $D_2-z$ and so its degree is the same in $D_1$ and $D$. 
By Lemma~\ref{lemma:2.4}, in $D_1$ there exists either an Ore-collapsible subdigraph $O$ or an emerald $E$ disjoint from $x$. 
In the former case $O$ is an Ore-collapsible subdigraph of $D$ disjoint from $T$, and in the later one $E$ is an emerald in $D$ 
disjoint from $T$ even if $y\in V(E)$ because $y$ has the same degree in $D_1$ and $D$. 
\end{proof}

\begin{restatable}{lemma}{twosix}
\label{lemma:2.6}
    If $R$ is an Ore-collapsible induced subdigraph of a $4$-Ore digraph $D$, then there exists a diamond or an
emerald of $D$ whose vertices lie in $V(R)$.
\end{restatable}
\begin{proof}
Let $D$ be a digraph.
Let $R$ be a minimal counterexample to this lemma, and let $\partial(R) = \{u,v\}$ and $H = D\langle R\rangle \cup [u,v]$. 
If $H=\bid{K_4}$, then $R$ is a diamond in $D$.
Suppose now that $H$ is the Ore-composition of two $4$-Ore digraphs $H_1$ (the digon side with replaced digon $[x,y]$) and $H_2$ (the split side with split vertex $z$).
If $\{u,v\} \not\subset V(H_2)$, then by Lemma~\ref{lemma:2.4} there exists an Ore-collapsible subdigraph in $H_2$ disjoint from $z$. As it is smaller
than $H$, it contains an emerald or a diamond as desired, a contradiction.

Now assume that $\{u,v\} \subset V(H_2)$, then $H_1$ is an Ore-collapsible subdigraph of $D$ smaller than $H$,
and by induction, $H_1$ contains a diamond or an emerald in $D$.
\end{proof}

\begin{lemma}\label{lemma:4-ore_diamond_emerald}
If $D$ is a $4$-Ore digraph and $v$ is a vertex in $D$, then $D$ contains
a diamond or an emerald disjoint from $v$.
\end{lemma}

\begin{proof}
Follows from Lemmas~\ref{lemma:2.4} and~\ref{lemma:2.6}.
\end{proof}

\begin{lemma}\label{lemma:4-ore_diamond_emerald_disjoint_triangle}
If $D$ is a $4$-Ore digraph and $T$ is a bidirected triangle in $D$,
then either $D=\bid{K_4}$ or $D$ contains a diamond or an emerald disjoint
from $T$.
\end{lemma}

\begin{proof}
Follows from Lemmas~\ref{lemma:2.5} and~\ref{lemma:2.6}.
\end{proof}

The following theorem was formulated for undirected graphs, but
by replacing every edge by a digon, it can be restated as follows:
\begin{theorem}[Kostochka and Yancey~\cite{kostochka2018brooks},~Theorem~6]\label{theorem:kostochka_yancey_4ore}
Let $D$ be a $4$-dicritical bidirected digraph. 

If $\frac{10}{3}n(D)-m(D) > 1$, then $D$ is $4$-Ore and $\frac{10}{3}n(D)-m(D) = \frac{4}{3}$.
\end{theorem}

\begin{lemma}\label{lemma:pot_4ore}
If $D$ is a $4$-Ore digraph with $n$ vertices, then $\rho(D) \leq \frac{4}{3}+\epsilon n - \delta\frac{2(n-1)}{3}$.
\end{lemma}

\begin{proof}
Follows from Theorem~\ref{theorem:kostochka_yancey_4ore} and Lemma~\ref{lemma:matching_ore}.
\end{proof}

\begin{lemma}[Kostochka and Yancey~\cite{kostochka2018brooks}, Claim~16]\label{lemma:pot_subgraph_4Ore}
Let $D$ be a $4$-Ore digraph. If $R \subseteq D$ and $0 < n(R) < n(D)$,
then $\frac{10}{3}n(R) - m(R) \geq \frac{10}{3}$.
\end{lemma}

\begin{sloppypar}
\begin{restatable}{lemma}{oldfourcliqueinfourore}
\label{lemma:old_4_clique_in_4Ore}
    Let $D$ be a $4$-Ore digraph obtained from a copy $J$ of $\bid{K_4}$ by successive Ore-compositions with $4$-Ore digraphs, 
vertices and digons in $J$ being always on the digon side.
Let $[u,v]$ be a digon in $D\langle V(J) \rangle$.
For every $3$-dicolouring $\phi$ of $D \setminus [u,v]$, vertices in $V(J)$ receive distinct colours except $u$ and $v$.
\end{restatable}
\end{sloppypar}
\begin{proof}
We proceed by induction on $n(D)$, the result holding trivially when $D$ is $\bid{K_4}$.
Now assume that $D$ is the Ore-composition of $D_1$, the digon side containing $J$, and $D_2$, with $D_1$ and $D_2$ being $4$-Ore digraphs.
Let $[x,y] \subseteq A(D_1)$ be the replaced digon in this Ore-composition, and let $z \in V(D_2)$ be the split vertex.
Let $\phi$ be a $3$-dicolouring of $D\setminus [u,v]$. Then $\phi$ induces a $3$-dicolouring of $D\langle V(D_2-z) \cup \{x,y\} \rangle$.
Necessarily $\phi(x) \neq \phi(y)$, for otherwise $\phi_2$ defined by $\phi_2(w) = \phi(w)$ if $w\in V(D_2-z)$ and $\phi_2(z)= \phi(x)$ is a $3$-dicolouring of $D_2$,
contradicting the fact that $4$-Ore digraphs have dichromatic number $4$ by Lemma~\ref{lemma:4ore_4_dicritical}.
Hence $\phi$ induces a $3$-dicolouring of $D_1 \setminus [u,v]$. So, by the induction hypothesis,
vertices in $V(J)$ have distinct colours in $\phi$, except $u$ and $v$.
\end{proof}

\begin{sloppypar}
\begin{restatable}{lemma}{oldfourcliqueinfourorebis}
\label{lemma:old_4_clique_in_4Ore_bis}
Let $D$ be a $4$-Ore digraph obtained from a copy $J$ of $\bid{K_4}$ by successive Ore-compositions with $4$-Ore digraphs, 
vertices and digons in $J$ being always on the digon side.
Let $v$ be a vertex in $V(J)$. For every $3$-dicolouring $\phi$ of $D - v$, vertices in $J$ receive distinct colours.
\end{restatable}
\end{sloppypar}
\begin{proof}
We proceed by induction on $n(D)$, the result holding trivially when $D$ is $\bid{K_4}$.
Now assume that $D$ is the Ore-composition of $D_1$, the digon side containing $J$, and $D_2$, with $D_1$ and $D_2$ being $4$-Ore digraphs.
Let $[x,y] \subseteq A(D_1)$ be the replaced digon in this Ore-composition, and let $z \in V(D_2)$ be the split vertex.
Let $\phi$ be a $3$-dicolouring of $D - v$. 
If $v \in \{x,y\}$, then $\phi$ is a $3$-dicolouring of $D_1- v$ and the result follows by induction. Now assume $v \not\in\{x,y\}$.
Then $\phi$ induces a $3$-dicolouring of $D\langle V(D_2-z) \cup \{x,y\} \rangle$.
Necessarily $\phi(x) \neq \phi(y)$, for otherwise $\phi_2$ defined by $\phi_2(w) = \phi(w)$ if $w\in V(D_2-z)$ and $\phi_2(z)= \phi(x)$ is a $3$-dicolouring of $D_2$,
contradicting the fact that $4$-Ore digraphs have dichromatic number $4$ by Lemma~\ref{lemma:4ore_4_dicritical}.
Hence $\phi$ induces a $3$-dicolouring of $D_1 - v$. So, by the induction hypothesis, vertices in $V(J)$ have distinct colours in $\phi$.
\end{proof}

\section{Proof of Theorem~\ref{thm:main_thm_potential}}\label{sec:main_proof}

Let $D$ be a $4$-dicritical digraph, $R$ be an induced subdigraph of $D$ with $4 \leq n(R) < n(D)$ and $\phi$ a $3$-dicolouring of $R$.
The {\bf $\phi$-identification} of $R$ in $D$, denoted by $D_{\phi}(R)$ is the digraph obtained from $D$ by identifying for each $i\in [3]$ 
the vertices coloured $i$ in $V(R)$ to a vertex $x_i$, adding the digons $[x_i,x_j]$ for all $1\leq i < j\leq 3$.
Observe that $D_{\phi}(R)$ is not $3$-dicolourable. Indeed, assume for a contradiction that $D_{\phi}(R)$ has a $3$-dicolouring $\phi'$. Since $\{x_1,x_2,x_3\}$ induces a $\bid{K_3}$, we may assume without loss of generality that $\phi'(x_i) =i$ for $i\in [3]$.
Consider the $3$-colouring $\phi''$ of $D$ defined by $\phi''(v)=\phi'(v)$ if $v \not\in R$ and $\phi''(v)=\phi(v)$ if $v \in R$. 
One easily checks that $\phi''$ is a $3$-dicolouring of $D$, a contradiction to the fact that $\vec{\chi}(D) = 4$.

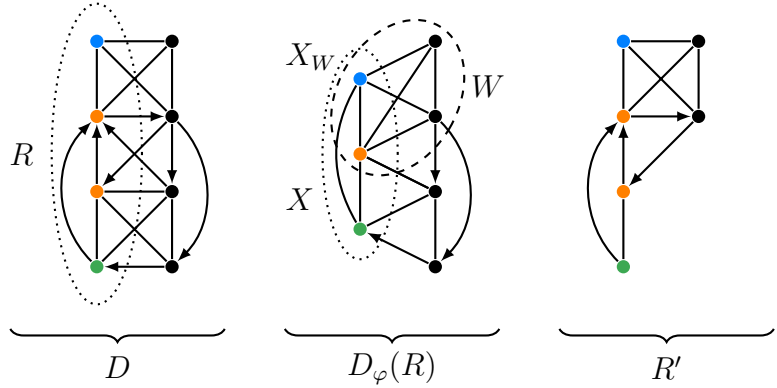
\begin{figure}[hbtp]
    \begin{minipage}{\linewidth}
        \begin{center}	
            \begin{tikzpicture}[thick,scale=1, every node/.style={transform shape}]
                \tikzset{vertex/.style = {circle,fill=black,minimum size=5pt,
                                    inner sep=0pt}}
                \tikzset{edge/.style = {->,> = latex}}
                \node[vertex,g-blue] (1l) at  (-1,0) {};
                \node[vertex,orange] (2l) at  (-1,-1) {};
                \node[vertex,orange] (3l) at  (-1,-2) {};
                \node[vertex,g-green] (4l) at  (-1,-3) {};
                \node[vertex] (1r) at  (0,0) {};
                \node[vertex] (2r) at  (0,-1) {};
                \node[vertex] (3r) at  (0,-2) {};
                \node[vertex] (4r) at  (0,-3) {};
                \draw[] (1l) -- (1r) -- (2r) -- (1l) -- (2l) -- (1r);
                \draw[edge] (2l) to (2r){};
                \draw[edge] (2r) to (3r){};
                \draw[edge] (3r) to (2l){};
                \draw[edge] (2r) to (3l){};
                \draw[edge] (3l) to (2l){};
                \draw[] (3l) -- (4r) -- (3r) -- (4l) -- (3l) --(3r);
                \draw[edge] (4r) to (4l){};
                \draw[edge] (4l) to[in=-135, out=135] (2l){};
                \draw[edge] (2r) to[in=45, out=-45] (4r){};
                \draw[dotted] (-1,-1.5) ellipse (0.6 and 2);
                \node[] (R) at (-2,-1.5){$R$};
                \draw [decorate,decoration={brace,amplitude=5pt,mirror,raise=4ex}] (-2.15, -3.1) -- (0.7,-3.1) node[midway,yshift=-3em]{$D$};
                \node[vertex,g-blue] (1lb) at  (2.5,-0.5) {};
                \node[vertex,orange] (2lb) at  (2.5,-1.5) {};
                \node[vertex,g-green] (3lb) at  (2.5,-2.5) {};
                \node[vertex] (1rb) at  (3.5,0) {};
                \node[vertex] (2rb) at  (3.5,-1) {};
                \node[vertex] (3rb) at  (3.5,-2) {};
                \node[vertex] (4rb) at  (3.5,-3) {};
                \draw[] (1lb) -- (2lb) -- (1rb) -- (2rb) -- (1lb) -- (1rb);
                \draw[] (3rb) -- (2lb) -- (2rb);
                \draw[] (4rb) -- (3rb);
                \draw[] (3lb) to[in=-120, out=120] (1lb){};
                \draw[] (3lb) -- (2lb) -- (3rb) -- (3lb);
                \draw[edge] (2rb) to (3rb){};
                \draw[edge] (2rb) to[in=45, out=-45] (4rb){};
                \draw[edge] (4rb) to (3lb){};
                \draw[dotted] (2lb) ellipse (0.5 and 1.4);
                \draw[dashed, rotate around={-30:(3,-0.75)}] (3,-0.75) ellipse (0.8 and 1.1);
                \node[] (X) at (1.7,-2.1){$X$};
                \node[] (X) at (1.85,-0.2){$X_W$};
                \node[] (W) at (4.2,-0.6){$W$};
                \draw [decorate,decoration={brace,amplitude=5pt,mirror,raise=4ex}] (1.5, -3.1) -- (4.35,-3.1) node[midway,yshift=-3em]{$D_\phi(R)$};
                \node[vertex,g-blue] (1lc) at  (6,0) {};
                \node[vertex,orange] (2lc) at  (6,-1) {};
                \node[vertex,orange] (3lc) at  (6,-2) {};
                \node[vertex,g-green] (4lc) at  (6,-3) {};
                \node[vertex] (1rc) at  (7,0) {};
                \node[vertex] (2rc) at  (7,-1) {};
                \draw[] (1lc) -- (1rc) -- (2rc) -- (1lc) -- (2lc) -- (1rc);
                \draw[edge] (2lc) to (2rc){};
                \draw[edge] (2rc) to (3lc){};
                \draw[edge] (3lc) to (2lc){};
                \draw[] (4lc) -- (3lc);
                \draw[edge] (4lc) to[in=-135, out=135] (2lc){};
                \draw [decorate,decoration={brace,amplitude=5pt,mirror,raise=4ex}] (1.5+3.65, -3.1) -- (1.5+3.65+2.85,-3.1) node[midway,yshift=-3em]{$R'$};
              \end{tikzpicture}
          \caption{A $4$-dicritical digraph $D$ together with an induced subdigraph $R$ of $D$ and $\phi$ a $3$-dicolouring of $R$, the $\phi$-identification $D_\phi(R)$ of $R$ in $D$ and the dicritical extension $R'$ of $R$ with extender $W$ and core $X
          _W$. For clarity, the digons are represented by undirected edges. }
          \label{fig:dicritical_extension}
        \end{center}    
      \end{minipage}
    \end{figure}

Now let $W$ be a $4$-dicritical subdigraph of $D_{\phi}(R)$ and $X=\{x_1, x_2, x_3\}$. Then we say that $R' = D \langle (V (W ) \setminus X)\cup R \rangle$ is the {\bf dicritical extension} of $R$ with {\bf extender} $W$. We call $X_W=X\cap V(W)$ the {\bf core} of the extension.
 Note that $X_W$ is not empty, because $W$ is not a subdigraph of $D$. Thus $1\leq |X_W| \leq 3 $. See Figure~\ref{fig:dicritical_extension} for an example of a $\phi$-identification and a dicritical extension.

\medskip

Let $D$ be a counterexample to Theorem~\ref{thm:main_thm_potential} with minimum number of vertices.
By Lemma~\ref{lemma:pot_4ore}, $D$ is not $4$-Ore. Thus $\rho(D) > 1$.

\begin{restatable}{claim}{firstclaim}\label{claim:general_bound_smaller_4_dicritical}
    If $\Tilde{D}$ is a $4$-dicritical digraph with $n(\Tilde{D}) < n(D)$, then $\rho(\Tilde{D}) \leq \frac{4}{3} + 4\epsilon - 2\delta$.
\end{restatable}
\begin{proofclaim}
    If $\Tilde{D}$ is not $4$-Ore, then $\rho(\Tilde{D}) \leq 1$ by minimality of $D$. Thus $\rho(\Tilde{D}) \leq \frac{4}{3}+4\epsilon-2\delta$ because $4\epsilon - 2\delta \geq \frac{-1}{3}$.
    Otherwise, by Lemma~\ref{lemma:pot_4ore}, $\rho(\Tilde{D}) \leq \frac{4}{3} + \epsilon n(\Tilde{D}) - \delta \frac{2(n(\Tilde{D})-1)}{3} \leq \frac{4}{3} + 4\epsilon - 2 \delta$ because $\delta \geq \frac{3}{2} \epsilon$ and $n(\Tilde{D}) \geq 4$. 
\end{proofclaim}

\begin{restatable}{claim}{secondclaim}\label{claim:first_potential_computation}
Let $R$ be a subdigraph of $D$ with $4 \leq n(R) < n(D)$. If $R'$ is a dicritical extension of $R$ with extender $W$ and core $X_W$, then
\[
\rho(R') \leq \rho(W) + \rho(R) - \left (\rho(\bid{K_{|X_W|}}) + \delta \cdot T(\bid{K_{|X_W|}})\right) + 
\delta \cdot (T(W)-T(W - X_W))
\]
and in particular
\[
\rho(R') \leq \rho(W) + \rho(R) - \frac{10}{3}-\epsilon+\delta.
\]
\end{restatable}
\begin{proofclaim}
We have
\begin{itemize}
    \item $n(R') = n(W) - |X_W|+n(R)$,
    \item $m(R') \geq m(W)+m(R)-m(\bid{K_{|X_W|}})$,
    \item $T(R') \geq T(W - X_W) + T(R)$
\end{itemize}
and by summing these inequalities, we get the first result.

Now observe that $T(W)-T(W - X_W)\leq |X_W|$ by Lemma~\ref{lemma:computation_T_minus_one_vertex}, and that the maximum
of $- \left(\rho(\bid{K_{|X_W|}}) + \delta T(\bid{K_{|X_W|}})\right) + \delta |X_W|$ is reached when $|X_W|=1$, in which case it is equal to $-\frac{10}{3}-\epsilon+\delta$. The second inequality follows.
\end{proofclaim}

\begin{restatable}{claim}{thirdclaim}\label{claim:pot_general_bad_bound}
If $R$ is a subdigraph of $D$ with $4 \leq n(R) < n(D)$, then $\rho(R) \geq \rho(D) + 2 -3\epsilon+\delta > 3 - 3\epsilon + \delta$.
\end{restatable}
\begin{proofclaim}
We proceed by induction on $n-n(R)$. Let $R'$ be a dicritical extension of $R$ with extender $W$ and core $X_W$.
By Claim~\ref{claim:first_potential_computation}, we have
\[
\rho(R') \leq \rho(W) + \rho(R) - \frac{10}{3}-\epsilon+\delta.
\]
Either $V(R')=V(D)$ and so $\rho(R') \geq \rho(D)$ or $V(R')$ is a proper subset of $V(D)$ and, since $R$ is a proper subdigraph of $R'$, by induction $\rho(R') \geq \rho(D) + 2 - 3\epsilon + \delta \geq \rho(D)$.
In both cases, $\rho(R') \geq \rho(D)$. Now $W$ is smaller than $D$ so $\rho(W) \leq \frac{4}{3}+4\epsilon-2\delta$ by Claim~\ref{claim:general_bound_smaller_4_dicritical}.
Thus 
\[
\rho(D) \leq \rho(R') \leq \frac{4}{3}+4\epsilon-2\delta + \rho(R) -\frac{10}{3}-\epsilon+\delta.
\]
This gives $\rho(R) \geq \rho(D) + 2 -3\epsilon+\delta > 3 - 3\epsilon + \delta$, because $\rho(D) > 1$.
\end{proofclaim}

As a consequence of Claim~\ref{claim:pot_general_bad_bound}, any subdigraph (proper or not) of size at least $4$ has potential at least $\rho(D)$.

We say that an induced subdigraph $R$ of $D$ is {\bf collapsible} if, for every $3$-dicolouring $\phi$ of $R$, its dicritical extension $R'$ (with extender $W$ and core $X_W$) is $D$, has core of size $1$ (i.e. $|X_W|=1$), and the border $\partial_D(R)$ of $R$ is monochromatic in $\phi$.

\begin{restatable}{claim}{fourthclaim}\label{claim:pot_not_monochromatic_boundary}
    Let $R$ be an induced subdigraph of $D$ and $\phi$ a 3-dicolouring of $R$ such that $\partial (R)$ is not monochromatic in $\phi$. If $D$ is a dicritical extension of $R$ dicoloured by $\phi$ with extender $W$ and core $X_W$ with $|X_W|=1$, then
    \[
        \rho(R) \geq \rho(D) +3 - 3\epsilon + \delta.
    \] 
\end{restatable}
\begin{proofclaim}
Assume $D$ is a dicritical extension of $R$ dicoloured by $\phi$ with extender $W$ and core $X_W$ with $|X_W|=1$.
Observe that each of the following inequalities holds:
    \begin{itemize}
        \item $n(D) = n(W) - |X_W| + n(R) = n(W)+n(R) - 1$,
        \item $m(D) \geq m(W)+m(R)-m(\bid{K_{|X_W|}}) + 1 = m(W)+m(R)+1$ because $\partial_D(R)$ is not monochromatic in $\phi$, and
        \item $T(D) \geq T(W - X_W) + T(R) \geq T(W) + T(R) - 1$ by Lemma~\ref{lemma:computation_T_minus_one_vertex}.
    \end{itemize}
By Claim~\ref{claim:general_bound_smaller_4_dicritical}, we have
\[
\rho(D) \leq  \rho(W) + \rho(R) - \left(\frac{10}{3}+\epsilon\right) - 1 +  \delta \leq \left(\frac{4}{3}+4\epsilon - 2\delta\right) + \rho(R) - \frac{13}{3} -\epsilon + \delta
\]
and so $\rho(R) \geq \rho(D) + 3 - 3\epsilon + \delta$.
\end{proofclaim}

\begin{restatable}{claim}{fifthclaim}\label{claim:non_collapsible_pot_bound}
    If $R$ is a subdigraph of $D$ with $4 \leq n(R) < n(D)$ and $R$ is not collapsible, then
    $\rho(R) \geq \rho(D) + \frac{8}{3}-\epsilon-\delta > \frac{11}{3}-\epsilon-\delta$.
\end{restatable}
\begin{proofclaim}
Let $R'$ be a dicritical extension of $R$ dicoloured by $\phi$ with extender $W$ and core $X_W$.

\medskip

(i) If $R'$ is not $D$, then it has a dicritical extension $R''$ with extender $W'$. 
By (the consequence of) Claim~\ref{claim:pot_general_bad_bound}, we have 
$\rho(D) \leq \rho(R'')$ ; by Claim~\ref{claim:first_potential_computation} (applied twice), $\rho(R'') \leq \rho(R) + \rho(W') + \rho(W) + 2\left(-\frac{10}{3} - \epsilon + \delta\right)$ ; both $W$ and $W'$ are smaller than $D$, so, by Claim~\ref{claim:general_bound_smaller_4_dicritical}, $\rho(W), \rho(W') \leq \frac{4}{3}+4\epsilon-2\delta$.
Those three inequalities imply
\[
\rho(D) \leq \rho(R'') \leq \rho(R) + 2 \left( \frac{4}{3}+4\epsilon-2\delta \right) + 2\left(-\frac{10}{3} - \epsilon + \delta\right) = \rho(R)  -4  + 6\epsilon - 2\delta
\]

and so $\rho(R) \geq \rho(D) + 4 - 6\epsilon + 2\delta \geq \rho(D) + \frac{8}{3}-\epsilon-\delta$. 

\medskip
(ii) If $R'=D$ and $|X_W|=2$, then 
$\rho(\bid{K_{|X_W|}}) + \delta T(\bid{K_{|X_W|}}) =  \frac{14}{3} +2\epsilon$, and, by Lemma~\ref{lemma:computation_T_minus_one_vertex},  $T(W)-T(W - X_W) \leq |X_w| = 2$.
Thus, by Claim~\ref{claim:first_potential_computation},
\[ \rho(D) \leq  \rho(W) + \rho(R) - \frac{14}{3} -2\epsilon + 2\delta
\]
Now, since $W$ is smaller than $D$, $\rho(W) \leq \frac{4}{3}+4\epsilon-2\delta$ by Claim~\ref{claim:general_bound_smaller_4_dicritical}. Thus
\[
\rho(D) \leq \rho(R) + \frac{4}{3}+4\epsilon-2\delta - \frac{14}{3} -2\epsilon + 2\delta
= \rho(R) - \frac{10}{3} +2\epsilon 
\]
and so $\rho(R) \geq \rho(D) + \frac{10}{3} -2\epsilon \geq \rho(D) + \frac{8}{3}-\epsilon-\delta$.

\medskip

(iii) If $R'=D$ and $|X_W|=3$, then
$\rho(\bid{K_{|X_W|}}) + \delta T(\bid{K_{|X_W|}}) =  4 +3\epsilon$, and, by Lemma~\ref{lemma:computation_T_minus_one_vertex}, $T(W)-T(W - X_W) \leq |X_w| = 3$.
Thus, by Claim~\ref{claim:first_potential_computation},
\[
 \rho(D) \leq \rho(W) + \rho(R) - 4 - 3\epsilon + 3\delta.
\]
Now, since $W$ is smaller than $D$, $\rho(W) \leq \frac{4}{3}+4\epsilon-2\delta$ by Claim~\ref{claim:general_bound_smaller_4_dicritical}. Thus
\[\rho(D) \leq \rho(R) + \frac{4}{3}+4\epsilon-2\delta - 4 -3\epsilon + 3\delta = \rho(R)  - \frac{8}{3} + \epsilon + \delta
\]
and so $\rho(R) \geq \rho(D) + \frac{8}{3} -\epsilon - \delta$.

\medskip

(iv) If $R'=D$, $|X_W|=1$ and $\partial(R)$ is not monochromatic in $\phi$, then, by Claim~\ref{claim:pot_not_monochromatic_boundary}, we have $\rho(R) \geq \rho(D) + 3 - 3\epsilon + \delta \geq \rho(D) + \frac{8}{3} - \epsilon - \delta$.

\medskip
If $R$ is not collapsible, then, by definition, it has a dicritical extension $R'$ satisfying the hypothesis of one of the cases (i)--(iv). In any case, $\rho(R) \geq \rho(D) + \frac{8}{3} - \epsilon - \delta$.
\end{proofclaim}

Recall that a {\bf $k$-cutset} in a graph $G$ is a set $S$ of $k$ vertices such that $G-S$ is not connected. A graph is {\bf $k$-connected} if it has more than $k$ vertices and has no $(k-1)$-cutset.
A {\bf $k$-cutset} in a digraph is a $k$-cutset in its underlying graph, and a digraph is {\bf $k$-connected} if its underlying graph is $k$-connected.

\begin{restatable}{claim}{sixthclaim}\label{claim:2_connected}
    $D$ is $2$-connected.
\end{restatable}
\begin{proofclaim}
    Suppose for contradiction that $\{x\}$ is a $1$-cutset in $D$. Let $(A_0, B_0)$ be a partition of $V(D - x)$ into non-empty sets such that there is no edge between $A_0$ and $B_0$, and set $A=A_0\cup \{x\}$ and $B=B_0\cup \{x\}$.

    Since $D$ is 4-dicritical, there exist a $3$-dicolouring $\phi_A$ of $D\langle A\rangle$ and a $3$-dicolouring $\phi_B$ of $D\langle B\rangle$. Free to swap the colours, we may assume $\phi_A(x) = \phi_B(x)$. Let $\phi$ be defined by $\phi(v) = \phi_A(v)$ if $v\in A$ and $\phi(v) = \phi_B(v)$ if $v\in B$. Since $\dic(D)= 4$, $D$, coloured with $\phi$, must contain a monochromatic directed cycle. Such a directed cycle must be contained in $D\langle A\rangle$ or $D\langle B\rangle$, a contradiction.
\end{proofclaim}

\begin{restatable}{claim}{seventhclaim}\label{claim:3_connected}
    $D$ is $3$-connected. In particular, $D$ contains no diamond.
\end{restatable}
\begin{proofclaim}
Suppose for contradiction that $\{x,y\}$ is a $2$-cutset of $D$. Let $(A_0, B_0)$ be a partition of $V(D)\setminus \{x,y\}$ into non-empty sets such that there is no edge between $A_0$ and $B_0$, and set $A=A_0\cup \{x,y\}$ and $B=B_0\cup \{x,y\}$.

Assume for a contradiction that there exists a $3$-dicolouring $\phi_A$ of $D\langle A\rangle$ and a $3$-dicolouring $\phi_B$ of $D\langle B\rangle$  such that $\phi_A(x)  \neq \phi_A(y)$ and $\phi_B(x) \neq \phi_B(y)$.
Free to swap the colours, we may assume $\phi_A(x) = \phi_B(x)$ and $\phi_A(y)  = \phi_B(y)$.
Let $\phi$ be defined by $\phi(v) = \phi_A(v)$ if $v\in A$ and $\phi(v) = \phi_B(v)$ if $v\in B$. Every directed cycle either is in $D\langle A\rangle$, or is in $D\langle B\rangle$ or contains both $x$ and $y$. Therefore it cannot be monochromatic with $\phi$ because $\phi_A$ and $\phi_B$ are $3$-dicolourings of $D\langle A\rangle$ and $D\langle B\rangle$ respectively, and $\phi(x) \neq \phi(y)$. Thus 
$\phi$ is a $3$-dicolouring of $D$, a contradiction.
Henceforth either $D\langle A\rangle$ or $D\langle B\rangle$ has no $3$-dicolouring $\phi$ such that
$\phi(x) \neq \phi(y)$. Suppose without loss of generality that it is $D\langle A\rangle$.

Let $D_A = D\langle A\rangle \cup [x,y]$.
$D_A$ is not $3$-dicolourable because in every $3$-dicolouring of $D\langle A\rangle$,
$x$ and $y$ are coloured the same.
Let $D_B$ be the digraph obtained from $D\langle B\rangle$ by identifying $x$ and $y$ into a vertex $z$.  Assume for a contradiction that $D_B$ has a $3$-dicolouring $\psi_{B}$. Set $\psi(x)=\psi(y)=\psi_{B}(z)$, and $\psi(u)=\psi_{B}(u)$ for every $u \in B \setminus\{x,y\}$.
Then consider a $3$-dicolouring $\psi_A$ of
$D\langle A\rangle$ such that $\psi_A(x)=\psi(x)=\psi_A(y)=\psi(y)$ (such a colouring exists because $A$ is a proper subdigraph of $D$) and we set $\psi(u) = \psi_A(u)$ for every $u\in V(A) \setminus \{ x,y\}$.
As $D$ is not $3$-dicolourable, it contains a monochromatic directed cycle $C$ (with respect to $\psi$). The cycle $C$ is not included in $D\langle A\rangle$ nor in $D_B$. As a consequence, there is a monochromatic
directed path from $\{x,y\}$ to $\{x,y\}$ in $B$, and so there is a monochromatic directed
cycle in $D_B$ for $\psi_{B}$, a contradiction.
Therefore $D_B$ is not $3$-dicolourable 

Now $D_A$ has a $4$-dicritical subdigraph $W_A$ which necessarily contains $\{x,y\}$, and  $D_B$ has a $4$-dicritical subdigraph $W_B$ which necessarily contains $z$.
As $W_A$ and $W_B$ are $4$-dicritical digraphs smaller than $D$, we have $\rho(W_A),\rho(W_B) \leq \frac{4}{3}+4\epsilon-2\delta$ by Claim~\ref{claim:general_bound_smaller_4_dicritical}.
Let $H$ be the subdigraph of $D$ induced by $V(W_A)\cup V(W_B -z)$.

Note that $n(H) = n(W_A) + n(W_B) -1$ and $m(H) \geq m(W_A) + m(W_B) -2$.
Moreover $T(H) \geq T(W_A - x) + T(W_B - z) \geq T(W_A) + T(W_B) -2$, by Lemma~\ref{lemma:computation_T_minus_one_vertex}.
Hence we have

\begin{eqnarray}
\rho(H) 
&\leq & \rho(W_A) + \rho(W_B) - \left(\frac{10}{3}+\epsilon\right) + (m(W_A) + m(W_B) - m(H)) + 2\delta \nonumber \\
&\leq  &\rho(W_A) + \rho(W_B) - \frac{10}{3} -\epsilon  +2 + 2\delta \nonumber \\
&=  &\rho(W_A) + \rho(W_B) - \frac{4}{3} -\epsilon + 2\delta \label{eq:utile} \\
&\leq  &2\left (\frac{4}{3}+4\epsilon-2\delta\right) - \frac{4}{3} - \epsilon + 2\delta  \nonumber \\
& = &    \frac{4}{3} + 7\epsilon - 2\delta \nonumber
\end{eqnarray}

By Claim~\ref{claim:pot_general_bad_bound}, if $n(H) < n(D)$ then $\rho(H) > 3 - 3\epsilon + \delta$.
As $10\epsilon-3\delta \leq \frac{5}{3}$, we deduce that $H=D$. In the above chain of inequalities, we upper bounded $m(W_A) + m(W_B) - m(H)$ by $2$, doing the same computation with  $2+ m(W_A) + m(W_B) - m(H) -2$ instead of $2$, we get
$1 < \rho(D) = \rho(H) \leq \frac{4}{3} + 7\epsilon - 2\delta + (m(W_A) + m(W_B) - m(H)-2)$ and so $m(H) = m(W_A) + m(W_B) - 2$ because $2\delta -7\epsilon \leq \frac{2}{3}$.
In particular, there is no arc between $x$ and $y$ in $D$. Moreover, no arc was suppressed when identifying $x$ and $y$ into $z$ to obtain $D_B$, so $x$ and $y$ have no common out-neighbour (resp. in-neighbour) in $B_0$.

We first show that either $W_A$ or $W_B$ is not $4$-Ore.
Assume for contradiction that both $W_A$ and $W_B$ are $4$-Ore.
If $H=D$ is not bidirected, then by Lemma~\ref{lemma:ore_composition_only_interesting_if_bidirected}, $D\langle B\rangle$ admits a $3$-dicolouring $\phi_B$ such that $\phi_B(x) = \phi_B(y)$. Now let $\phi_A$ be a $3$-dicolouring of $D\langle A \rangle$. We have $\phi_A(x)=\phi_A(y)$. Free to exchange colours, we may assume
$\phi_A(x)=\phi_A(y) = \phi_B(x) = \phi_B(y)$.
Hence we can define the $3$-colouring $\phi$ of $D$ by $\phi(v)=\phi_A(v)$ if $v\in A$, and $\phi(v)=\phi_B(v)$ if $v\in B$. Observe that, since $A$ is bidirected, all neighbours of $x$ and $y$ in $D\langle A \rangle$ have a colour distinct from $\phi(x)$. Therefore there is no monochromatic directed cycle in $D$ coloured by $\phi$. Thus
$\phi$ is a $3$-dicolouring of $D$, a contradiction.
Therefore, $H=D$ is bidirected, and so $H$ is an Ore-composition of $W_A$ and $W_B$ (because $D$ is $2$-connected by Claim~\ref{claim:2_connected}), and so $D$ is $4$-Ore, a contradiction.
Henceforth, we may assume that either $W_A$ or $W_B$ is not $4$-Ore.

If none of $W_A$ and $W_B$ are a $4$-Ore, then by minimality of $D$, 
$\rho(W_A) \leq 1$ and $\rho(W_B) \leq 1$. Together with Equation~\eqref{eq:utile}, this yields 
\[
\rho(H) \leq \frac{2}{3} - \epsilon  + 2\delta \leq 1
\]
because $2\delta - \epsilon \leq \frac{1}{3}$, a contradiction.

If none of $W_A$ and $W_B$ is $\bid{K_4}$, then $\rho(W_A)+\rho(W_B) \leq 1 + (\frac{4}{3}+7\epsilon - 4\delta)$ (recall that if a digraph is $4$-Ore but not $\bid{K_4}$, then it has potential at most $\frac{4}{3}+7\epsilon-4\delta$ by Lemma~\ref{lemma:pot_4ore}). Thus, with Equation~\eqref{eq:utile}, we get  
\[
\rho(H) \leq 1 + \left(\frac{4}{3}+7\epsilon - 4\delta\right) - \frac{4}{3} - \epsilon  + 2\delta = 1 +6\epsilon -2\delta \leq 1
\]
because $\delta \geq 3\epsilon$. 

 Finally, if exactly one of $W_A$ or $W_B$ is isomorphic to $\bid{K_4}$, then either 
$T(W_A - x)  = T(W_A) = 2$ (if $W_A = \bid{K_4}$) or $T(W_B - z)  = T(W_B) = 2$ (if $W_B = \bid{K_4}$).
Therefore $T(H) \geq T(W_A - x) + T(W_B-z) \geq T(W_A) + T(W_B) - 1$ by Lemma~\ref{lemma:computation_T_minus_one_vertex}, and so 
\[
\rho(H) \leq  \rho(W_A)+\rho(W_B) - \left(\frac{10}{3}+\epsilon\right) + 2 + \delta.
\]
Now the non $4$-Ore digraph among $W_A$, $W_B$ has potential at most $1$ and the other has potential $\rho(\bid{K_4}) = \frac{4}{3} +4\epsilon - 2\delta$. Thus
\[
\rho(H) \leq 1 + \left(\frac{4}{3} +4\epsilon - 2\delta\right) - \left(\frac{10}{3}+\epsilon\right) + 2 + \delta = 1 + 3\epsilon - \delta \leq 1
\]
because $\delta \geq 3\epsilon$.

In all three cases, $\rho(D) = \rho(H) \leq 1$, which is a contradiction. Hence $D$ is $3$-connected.
\end{proofclaim}

\begin{restatable}{claim}{eighthclaim}\label{claim:4-Ore_collapsible}
If $R$ is a collapsible subdigraph of $D$, $u,v$ are in the boundary of $R$ and $D\langle R\rangle\cup[u,v]$ is $4$-Ore, then there exists $R' \subseteq R$ such that
\begin{enumerate}[label=(\roman*)]
    \item either $R'$ is an Ore-collapsible subdigraph of $D$, or
    \item $R'$ is an induced subdigraph of $R$, $n(R') < n(R)$, and there exist $u',v'$ in $\partial_D(R')$
    such that $R'\cup [u',v']$ is $4$-Ore.
\end{enumerate}
\end{restatable}
\begin{proofclaim}
If $\partial(R) = \{u,v\}$, then $R$ is Ore-collapsible and we are done.
Suppose now that there exists $w \in \partial(R)$ distinct from $u$ and $v$.
Let $H=D\langle R\rangle \cup [u,v]$. Observe that $H \neq \bid{K_4}$ as $u,v$ and $w$ receive
the same colour in any $3$-dicolouring of $D\langle R\rangle $ because $R$ is collapsible. Hence $H$ is the Ore-composition of two $4$-Ore digraphs
$H_1$ (the digon side with replaced digon $[x,y]$) and $H_2$ (the split side with split vertex $z$).

If $u$ or $v$ is in $V(H_2)$, then $R' = D\langle V(H_1)\rangle$ with $u'=x,v'=y$ satisfies (ii).
Now we assume that $u,v \in V(H_1) \setminus V(H_2)$.
By repeating this argument successively on $H_1$, and then on the digon-side of
$H_1$, etc, either we find a subdigraph $R'$ satisfying (ii) or $u$ and $v$ are in a copy $J$ of $\bid{K_4}$ such that $H$ is obtained by Ore-compositions between $J$ and some $4$-Ore digraphs with $J$ being always in the digon side.

Observe that $w \not\in V(J)$ because in any $3$-dicolouring of $H \setminus [u,v]$,
vertices in $J$ receive different colours by Lemma~\ref{lemma:old_4_clique_in_4Ore}, except $u$ and $v$. 
Hence at one step in the succession of Ore-compositions, $w$ was in the split-side
$S$ when a digon $e$ in $J$ has been replaced.
However $e \neq [u,v]$, so either $u$ or $v$ is not in $e$. Suppose without loss
of generality that $e$ is not incident to $v$.

We claim that $H'=R - v \cup [u,w]$ is not $3$-dicolourable. Otherwise, let $\phi$ be a $3$-dicolouring of $H'$. Then $\phi$ is a $3$-dicolouring of $H-v$ with $H$ $4$-Ore, so vertices in $J-v$ must receive pairwise different colours by Lemma~\ref{lemma:old_4_clique_in_4Ore_bis}. 
Let $\phi'$ be a $3$-dicolouring of $R$. Without loss of generality, we may assume that $\phi(x)=\phi'(x)$ for every $x \in V(J-v)$. If $y \in S$, let $\phi''(y)=\phi(y)$, and let $\phi''(y)=\phi'(y)$ if $y \not\in S$. Then $\phi''$ is a $3$-dicolouring of $R$ but with $\phi(u) \neq \phi(w)$, contradicting the fact that $R$ is collapsible. This shows that $H'=R - v \cup [u,w]$ is not $3$-dicolourable.

Hence $R-v\cup[u,w]$ contains a $4$-dicritical digraph $K$. By Lemma~\ref{lemma:pot_subgraph_4Ore}, $R' = D\langle V(K)\rangle$, as a subdigraph of $H$ which is a 4-Ore, satisfies $\frac{10}{3}n(R') - m(R') \geq \frac{10}{3}$. This implies that $\frac{10}{3}n(K)-m(K) \geq \frac{4}{3}$. Note also that $K$ is bidirected because $R-v$ is bidirected.
Thus, by Theorem~\ref{theorem:kostochka_yancey_4ore}, $K$ is $4$-Ore. Hence $R'$ with $u,w$ satisfies (ii).
\end{proofclaim}

\begin{restatable}{claim}{ninethclaim}\label{claim:no_collapsible}
    If $R$ is a subdigraph of $D$ with $n(R) < n(D)$ and $u,v \in V(R)$, then $R \cup [u,v]$ is $3$-dicolourable.
    As a consequence, there is no collapsible subdigraph in $D$.
\end{restatable}
\begin{proofclaim}
Assume for a contradiction that the statement is false. Consider a smallest induced subdigraph $R$ for which the statement does not hold. Then $K=R \cup [u,v]$ is 4-vertex-dicritical, that is for every vertex $v \in V(K), \dic(K-v) < 4=\dic(K)$. Note that $4$-vertex-dicritical digraphs smaller than $D$ satisfy the outcome of Theorem~\ref{thm:main_thm_potential} since adding arcs does not increase the potential. Note that $\rho(R) \leq \rho(K) + 2 + \delta$.

If $R$ is not collapsible, then, by Claim~\ref{claim:non_collapsible_pot_bound},
$\rho(R) \geq \rho(D) + \frac{8}{3}-\epsilon - \delta > \frac{11}{3} -\epsilon - \delta$.
But we also have $\rho(R) \leq \rho(K) + 2 + \delta \leq \frac{10}{3} + 4\epsilon - \delta$ by Claim~\ref{claim:general_bound_smaller_4_dicritical}, 
which is a contradiction because 
$5\epsilon \leq \frac{1}{3}$.
Hence $R$ is collapsible. 

 Let $\phi$ be a $3$-dicolouring of $R$. Observe that $\phi(u) = \phi(v)$ for otherwise $R \cup [u,v]$ would be $3$-dicolourable. 
Let $R'$ be the dicritical extension of $R$ with extender $W$ and core $X_W$. We have $R'=D$ and $|X_W|=1$. Since $R$ is collapsible, for every two vertices $u',v'$ on the boundary of $R$, $R \cup [u',v']$ is not $3$-dicolourable. Hence, free to consider $u',v'$ instead of $u,v$, we can suppose that $u$ and $v$ are on the boundary of $R$.
If $K$ is $4$-Ore, then, by Claim~\ref{claim:4-Ore_collapsible} and by minimality
of $R$, we have that $R$ is Ore-collapsible, and so has boundary of size $2$.
This contradicts the fact that $D$ is $3$-connected.
Hence $K$ is not $4$-Ore.

By Claim~\ref{claim:first_potential_computation}, we have
\begin{eqnarray*}
1 < \rho(D) = \rho(R')  & \leq &  \rho(W) + \rho(R) - \frac{10}{3} - \epsilon + \delta\\
& \leq & \rho(W) + (\rho(K)+2 + \delta) - \frac{10}{3} - \epsilon + \delta
\end{eqnarray*}
and as 
$\rho(K) \leq 1$
(because it is not $4$-Ore and by minimality of $D$) we get
\[
1 < 1 + \rho(W) - \left(\frac{4}{3}+ \epsilon-2\delta\right)
\]
that is $\rho(W) > \frac{4}{3}+ \epsilon-2\delta$.
But as $W$ is smaller than $D$, it satisfies Theorem~\ref{thm:main_thm_potential}. Thus, since 
$\epsilon -2\delta \geq \frac{-1}{3}$
, $W$ must be 4-Ore. Moreover, $W$ must be isomorphic to $\bid{K_4}$, for otherwise $\rho(W)$ would be at most $\frac{4}{3}+7\epsilon -4\delta$, and $\frac{4}{3}+7\epsilon -4\delta \geq \rho(W) > \frac{4}{3} + \epsilon - 2\delta$ would contradict $\delta \geq 3\epsilon$.
Hence $\rho(W)= \rho(\bid{K_4}) =\frac{4}{3} + 4\epsilon - 2\delta$ and $T(W - X_W) = 2 = T(W)$. Thus, by Claim~\ref{claim:first_potential_computation} and because $\delta \geq 3\epsilon$, we have 
\[
1 < \rho(D) \leq \rho(W) + \rho(K) + 2 + \delta - \frac{10}{3} - \epsilon \leq \rho(K) + 3\epsilon - \delta \leq \rho(K) \leq 1,
\]
 a contradiction.

This implies that $D$ does not contain any collapsible subdigraph. Indeed, assume for a contradiction that $D$ contains a collapsible subdigraph $R$, and let $u,v$ be two vertices in its boundary. Then there exists a 3-dicolouring $\phi$ of $R\cup [u,v]$, for which $\partial (R)$ is not monochromatic, a contradiction.
\end{proofclaim}

\begin{restatable}{claim}{tenthclaim}\label{claim:minus_one_vertex_plus_two_arcs}
    If $R$ is a subdigraph of $D$ with $n(R) < n(D)$ and $u,v,u',v' \in R$, then
$R \cup \{uv, u'v'\}$ is $3$-dicolourable. In particular, $D$ contains no copy of $\bid{K_4}$ minus two arcs.
\end{restatable}
\begin{proofclaim}
Assume for a contradiction that the statement is false. Consider a smallest subdigraph $R$ for which the statement does not hold. Then $K = R \cup \{uv, u'v'\}$ is $4$-dicritical and smaller than $D$, so $\rho(K) \leq \frac{4}{3} + 4\epsilon-2\delta$ by Claim~\ref{claim:general_bound_smaller_4_dicritical}.
By Claim~\ref{claim:no_collapsible}, $R$ is not collapsible, so, by Claim~\ref{claim:non_collapsible_pot_bound}, we have
$\rho(R) \geq \rho(D) + \frac{8}{3} - \epsilon - \delta > \frac{11}{3} - \epsilon - \delta$.
But $\rho(R) \leq \rho(K) + 2 + 2\delta \leq \frac{10}{3} + 4\epsilon$,
which is a contradiction as 
$5\epsilon+\delta \leq \frac{1}{3}$.
\end{proofclaim}

For any $v \in V(D)$, we denote by $n(v)$ its number of neighbours, that is $n(v) = |N^+(u) \cup N^-(v)|$, and by $d(v)$ its number of incident arcs, 
that is $d(v) = d^+(v) + d^-(v)$.

\begin{restatable}{claim}{eleventhclaim}\label{claim:degree6_3_or_6_neighbours}
Vertices of degree $6$ in $D$ have either three or six neighbours.
\end{restatable}

\begin{proofclaim}
Let $x$ be a vertex of degree $6$.

If $n(x)=4$, then let $a,b,c,d$ be its neighbours such that
$N^+(x)=\{a,b,c\}$ and $N^-(x)=\{a,b,d\}$.
Consider $D' = D - x \cup dc$.
By Claim~\ref{claim:minus_one_vertex_plus_two_arcs}, $D'$ has a $3$-dicolouring
$\phi$.
If $|\phi(N^-(x))| < 3$, then choosing $\phi(x)$ in $\{1,2,3\} \setminus \phi(N^-(x))$,
we obtain a $3$-dicolouring of $D$, a contradiction.
Hence $\phi(N^-(x)) = \{1,2,3\}$. We set $\phi(x)=\phi(d)$.
As $D$ is not $3$-dicolourable, $D$ contains a monochromatic directed cycle $C$.
This cycle $C$ must contain the arc $dx$, and an out-neighbour $z$ of $x$. Since $\phi(a)$, $\phi(b)$ and $\phi(d)$ are all distinct, necessarily $z=c$.
But then $C-x\cup dc$ is a monochromatic directed cycle in $D'$, a contradiction.

Similarly, if $n(v)=5$, let $N^+(x)=\{a,b,c\}$ and $N^-(x)=\{a,d,e\}$,
and consider $D'=D-x\cup db \cup dc$. By Claim~\ref{claim:minus_one_vertex_plus_two_arcs},
$D'$ has a $3$-dicolouring $\phi$.
If $|\phi(N^-(x))| < 3$, then choosing $\phi(x)$ in $\{1,2,3\} \setminus \phi(N^-(x))$, we obtain a $3$-dicolouring of $D$, a contradiction. Hence $\phi(N^-(x)) = \{1,2,3\}$. We set $\phi(x)=\phi(d)$.
As $D$ is not $3$-dicolourable, there is a monochromatic directed cycle $C$, which
must contain the arc $dx$ and an out-neighbour $z$ of $x$. Note that $z$ must be $b$ or $c$ because $\phi(a)\neq \phi(d)$.
Then $C-x \cup dz$ is a monochromatic directed cycle in $D'$, a contradiction.
\end{proofclaim}

\begin{restatable}{claim}{twelvethclaim}\label{claim:triangle_two_vertices_degree_6_dont_exist}
There is no bidirected triangle containing two vertices of degree $6$.
In particular, $D$ contains no emerald.
\end{restatable}
\begin{proofclaim}
Suppose that $D\langle\{x,y,z\}\rangle=\bid{K_3}$ and $d(x)=d(y)=6$.
By Claim~\ref{claim:degree6_3_or_6_neighbours}, $x$ and $y$ have exactly three
neighbours, and $N[x] \neq N[y]$ because $D$ contains no copy of $\bid{K_4}$ minus two arcs by Claim~\ref{claim:minus_one_vertex_plus_two_arcs}.
Let $u$ (resp. $v$) be the unique neighbour of $x$ distinct from $y$ and $z$
(resp. $x$ and $z$).
Consider $D' = D-\{x,y\}\cup[u,v]$. By Claim~\ref{claim:no_collapsible},
$D'$ has a $3$-dicolouring $\phi$.
Without loss of generality, suppose that $\phi(u)=1$ and $\phi(v)=2$.
If $\phi(z)=1$ (resp. $\phi(z)=2$, $\phi(z)=3$), we set $\phi(x)=2$ and $\phi(y)=3$ (resp. $\phi(x)=3$ and $\phi(y)=1$, $\phi(x)=2$ and $\phi(y)=1$).
In each case, this yields a $3$-dicolouring of $D$, a contradiction.
\end{proofclaim}

So now we know that $D$ contains no emerald, and no diamond by Claim~\ref{claim:3_connected}.

\begin{restatable}{claim}{thirteenthclaim}\label{claim:pot_subgraph_final}
If $R$ is an induced subdigraph of $D$ with $4 \leq n(R) < n(D)$, then $\rho(R) \geq \rho(D) + 3 +3\epsilon-3\delta$, except if $D-R$ contains a single vertex which has degree $6$ in $D$.
\end{restatable}

\begin{proofclaim}
Let $R$ be an induced subdigraph of $D$ with $4 \leq n(R) < n(D)$.
By Claim~\ref{claim:no_collapsible}, $R$ is not collapsible. 
Let $\phi$ be a 3-dicolouring of $R$, $R'$ be a dicritical extension of $R$ with extender $W$ and core $X_W$ (with respect to $\phi$). By (the consequence of) Claim~\ref{claim:pot_general_bad_bound},
we know that $\rho(R') \geq \rho(D)$.

\medskip

Assume first that $R' \neq D$. Then, by Claims~\ref{claim:pot_general_bad_bound} and~\ref{claim:first_potential_computation},
\[
\rho(D) + 2 - 3\epsilon+\delta \leq \rho(R') \leq \rho(W)+\rho(R) - \frac{10}{3} - \epsilon + \delta.
\]
Since $\rho(W) \leq \frac{4}{3}+4\epsilon-2\delta$ by Claim~\ref{claim:general_bound_smaller_4_dicritical}, we have $\rho(R) \geq \rho(D) + 4 -6\epsilon + 2 \delta \geq \rho(D) + 3 + 3\epsilon - 3\delta$, because $1 \geq 9\epsilon -5\delta$.
In the following we suppose that $R'=D$.
We distinguish three cases depending on the cardinality of $|X_W|$.
\medskip

\begin{itemize}
\item Assume first that $|X_W|=2$. Then, by Claim~\ref{claim:first_potential_computation} and Lemma~\ref{lemma:computation_T_minus_one_vertex},
\[
\rho(D) \leq \rho(R') \leq \rho(W) + \rho(R) - \frac{20}{3}-2\epsilon +2 +2\delta
\]
and, as $\rho(W) \leq \frac{4}{3}+4\epsilon-2\delta$ by Claim~\ref{claim:general_bound_smaller_4_dicritical}, 
we have $\rho(R) \geq \rho(D) + \frac{10}{3}-2\epsilon \geq \rho(D) + 3+3\epsilon-3\delta$ because $5\epsilon-3\delta \leq \frac{1}{3}$.

\medskip

\item Assume now that $|X_W|=3$.
If there is a vertex $v \in V(D-R)$ with two out-neighbours (resp. two in-neighbours) in $V(R)$ with the same colour for $\phi$, then 
\begin{itemize}
    \item $n(R') = n(W) - |X_W|+n(R)$,
    \item $m(R') \geq m(W)+m(R)-m(\bid{K_{|X_W|}}) + 1$ because $v$ has two in- or out-neighbour in $V(R)$ with the same colour for $\phi$,
    \item $T(R') \geq T(W - X_W) + T(R)$.
\end{itemize}
It follows that
\[
\rho(D) \leq \rho(R') \leq \rho(W)+\rho(R) - (10+3\epsilon-6)+3\delta - 1
\]
and so 
$\rho(R) \geq \rho(D) - \frac{4}{3} - 4\epsilon + 2\delta + 5 + 3\epsilon - 3 \delta \geq \rho(D) + \frac{11}{3}-\epsilon-\delta \geq \rho(D)+3+3\epsilon-3\delta$
because $4\epsilon - 2\delta \leq \frac{2}{3}$.
Now we assume that there is no vertex with two out-neighbours (resp. two in-neighbours) in $R$ with the same colour for $\phi$.
In other words, the in-degrees and out-degrees of vertices in $D-R$ are the same in $D$ and in $W$.

If $W$ is not $4$-Ore, then by Claim~\ref{claim:first_potential_computation}
\[
\rho(D) \leq \rho(R') \leq \rho(W) + \rho(R) - (10+3\epsilon-6) + 3\delta
\]
and, as $\rho(W) \leq 1$, we have $\rho(R) \geq \rho(D) + 3 +3\epsilon-3\delta$.

Now suppose $W$ is $4$-Ore.
If $W \neq \bid{K_4}$, then, by Lemma~\ref{lemma:4-ore_diamond_emerald_disjoint_triangle},
$W$ contains a diamond or an emerald disjoint from $X$, 
and this gives a diamond or an emerald in $D$ because
the degrees of vertices in $D-R$ are the same in $D$ and in $W$, which is a contradiction.

Now suppose that $W = \bid{K_4}$. Then $D - R$ has a single vertex of degree $6$ in $D$.

\item Assume finally that $|X_W|=1$. Since $R$ is not collapsible by Claim~\ref{claim:no_collapsible}, $\phi$ may have been chosen so that $\partial (R)$ is not monochromatic in $\phi$. Then, by Claim~\ref{claim:pot_not_monochromatic_boundary}, 
$\rho(R) \geq \rho(D) + 3 - 3\epsilon + \delta \geq \rho(D) + 3 + 3\epsilon - 3 \delta$, because $6\epsilon -4\delta \leq 0$.
\end{itemize}
\end{proofclaim}

In $D$, we say that a vertex $v$ is a {\bf simple in-neighbour} (resp. {\bf simple out-neighbour}) if $v$ is a in-neighbour (resp. out-neighbour) of $u$ 
and $[u,v]$ is not a digon in $D$. If $v$ is a simple in-neighbour or simple out-neighbour of $u$, we simply say that $v$ is a {\bf simple neighbour} of $u$.

\begin{restatable}{claim}{fourteenthclaim}\label{claim:degree7_7neighbours}
Vertices of degree $7$ have seven neighbours. In other words, every vertex of degree $7$ has only simple neighbours.
\end{restatable}

\begin{proofclaim}
Let $x$ be a vertex of degree $7$.
We suppose, without loss of generality, that $d^-(x)=3$ and $d^+(x)=4$.

\medskip

If $n(x) = 4$, then $x$ has a unique simple out-neighbour $a$.
As $D$ is $4$-dicritical, $D\setminus xa$ has a $3$-dicolouring $\phi$.
But then every directed cycle is either in $D\setminus xa$ or it contains $xa$ and thus an in-neighbour $t$ of $x$. In the first case, it is not monochromatic because $\phi$ is a $3$-dicolouring of $D\setminus xa$, and in the second case, it is not monochromatic because $[t,x]$ is a digon and so $\phi(t)\neq \phi(x)$.
Hence $\phi$ is a $3$-dicolouring of $D$, a contradiction.

\medskip

If $n(x)=5$, let $N^-(x)=\{a,b,c\}$ and $N^+(x)=\{a,b,d,e\}$.
By Claim~\ref{claim:minus_one_vertex_plus_two_arcs}, $D'=D-x \cup\{cd,ce\}$
has a $3$-dicolouring $\phi$.
If $|\phi(N^-(x))| < 3$, then choosing $\phi(x)$ in $\{1,2,3\} \setminus \phi(N^-(x))$ gives a $3$-dicolouring of $D$, a contradiction.
If $|\phi(N^-(x))| = 3$, then we set $\phi(x) = \phi(c)$. Suppose for a contradiction that there is a monochromatic directed cycle $C$ in $D$ (with $\phi$). Necessarily $C$ contains $x$
(since $\phi$ is a $3$-dicolouring of $D-x$) and so it must contain $c$ and one vertex $y$ in $\{d,e\}$ because $\phi(a)$, $\phi(b)$, and $\phi(c)$ are all distinct.
Then $C-x\cup cy$ is a monochromatic directed cycle in $D'$, a contradiction.
Therefore $\phi$ is a $3$-dicolouring of $D$, a contradiction.

\medskip

If $n(x)=6$, let $N^-(x)=\{a,b,c\}$ and $N^+(x)=\{a,d,e,f\}$.
Consider $D' = D - x \cup\{bd,be,bf\}$.

We first show that $D'$ is not $3$-dicolourable. Assume for a contradiction that there is a $3$-dicolouring $\phi$ of $D'$.
If $|\phi(N^-(x))| < 3$, then choosing $\phi(x)$ in $\{1,2,3\} \setminus \phi(N^-(x))$ gives a $3$-dicolouring of $D$, a contradiction.
Hence $|\phi(N^-(x))| = 3$. We set $\phi(x) = \phi(b)$. 
Since $D$ is not $3$-dicolourable, there exists a monochromatic directed cycle $C$ in $D$ (with $\phi$). Necessarily $C$ contains $x$ (since $\phi$ is a $3$-dicolouring of $D-x$) and so it must contain $b$ and one vertex $y$ in $\{d,e,f\}$ because $\phi(a)$, $\phi(b)$, and $\phi(c)$ are all distinct.
Then $C-x\cup by$ is a monochromatic directed cycle in $D'$, a contradiction.
This gives a $3$-dicolouring of $D$, a contradiction.

Henceforth $D'$ is not $3$-dicolourable, and so it contains a $4$-dicritical
digraph $\Tilde{D}$, smaller than $D$.
If $\Tilde{D}$ does not contain the three arcs $bd, be, bf$, then it can be obtained
from a proper induced subdigraph of $D$ by adding at most two arcs, and so it is
$3$-dicolourable by Claim~\ref{claim:minus_one_vertex_plus_two_arcs}, a contradiction.

Hence $\{b,d,e,f\}\subseteq V(\Tilde{D})$.
Now consider $U = D\langle V(\Tilde{D}) \cup \{x\}\rangle$.

\noindent $\bullet$ Assume first that $a \not\in V(U)$ or $c \not\in V(U)$. Then we have
\begin{itemize}
    \item $n(U) = n(\tilde{D}) +1$,
    \item $m(U) \geq m(\tilde{D}) + 1$ and 
    \item $T(U) \geq T(\tilde{D}-b)\geq T(\tilde{D})-1$ by Lemma~\ref{lemma:computation_T_minus_one_vertex}.
\end{itemize}
Hence 
\[
{
\everymath={\displaystyle}
\renewcommand{\arraystretch}{2}
\begin{array}{r l l}
    \rho(U) &\leq \rho(\Tilde{D}) + \frac{10}{3}+\epsilon-1+\delta                   & \\
            &\leq \frac{4}{3}+4\epsilon-2\delta + \frac{10}{3}+\epsilon - 1 + \delta &\text{by Claim~\ref{claim:general_bound_smaller_4_dicritical},} \\
            &=    1 + \frac{8}{3} +5\epsilon - \delta   & \\
            &<    \rho(D) + \frac{8}{3} +5\epsilon - \delta                      & \\
            &\leq \rho(D) + 3 + 3\epsilon-3\delta                      &\text{because $\frac{1}{3}\geq 2\delta+2\epsilon$}.
\end{array}
}
\]
Hence by Claim~\ref{claim:pot_subgraph_final}, $D-U$ has a single vertex of degree $6$ (in $D$), which must be either $a$ or $c$. Then we have
\begin{itemize}
    \item $n(D) = n(\Tilde{D}) +  2$,
    \item $m(D) \geq m(\Tilde{D}) - 3 + 11$ and
    \item $T(D) \geq T(\Tilde{D}-b) \geq T(\Tilde{D})-1$.
\end{itemize}
Thus 
\[
{
\everymath={\displaystyle}
\renewcommand{\arraystretch}{2}
\begin{array}{r l l}
    \rho(D) &\leq \rho(\Tilde{D}) + 2\left(\frac{10}{3}+\epsilon\right) - 8 + \delta            & \\
            &\leq \left(\frac{4}{3}+4\epsilon-2\delta\right) -\frac{4}{3} + 2\epsilon + \delta     & \text{by Claim~\ref{claim:general_bound_smaller_4_dicritical},} \\
            &\leq 1    &\text{because $6\epsilon -\delta \leq 1$.}
\end{array}
}
\]
This is a contradiction.

\noindent $\bullet$ Assume now that $a,c \in V(U)$, then we have 
\begin{itemize}
    \item $n(U) = n(\tilde{D}) +1$,
    \item $m(U) \geq m(\tilde{D}) +4$ and
    \item $T(U) \geq T(\tilde{D}-b)\geq T(\tilde{D})-1$ 
by Lemma~\ref{lemma:computation_T_minus_one_vertex}. 
\end{itemize}  
Thus 
\[
{
\everymath={\displaystyle}
\renewcommand{\arraystretch}{2}
\begin{array}{r l l}
    \rho(U) &\leq \rho(\Tilde{D}) + \frac{10}{3}+\epsilon - 4 +\delta           & \\
            &\leq \left(\frac{4}{3}+4\epsilon-2\delta\right) + \frac{10}{3}+\epsilon - 4 + \delta     & \text{ by Claim~\ref{claim:general_bound_smaller_4_dicritical},} \\
            &\leq 1 &\text{ because $5\epsilon - \delta \leq \frac{1}{3}$}.
\end{array}
}
\]
 Together with the consequence of Claim~\ref{claim:pot_general_bad_bound}, we get that $\rho(D) \leq \rho(U)  \leq  1$, a contradiction.
\end{proofclaim}

The {\bf $8^+$-valency} of a vertex $v$, denoted by $\nu(v)$, is the number of arcs incident to $v$ and a vertex of degree at least $8$. 

Let $D_6$ be the subdigraph of $D$ induced by the vertices of degree $6$ incident to digons.
Let us describe the connected components of $D_6$ and their neighbourhoods.
Remember that vertices of degree $7$ are incident to no digon by Claim~\ref{claim:degree7_7neighbours}, and so they do not have neighbours in $V(D_6)$.
If $v$ is a vertex in $D_6$, we define its {\bf neighbourhood valency} to be the sum of the $8^+$-valency of its neighbours of degree at least $8$.
We denote the neighbourhood valency of $v$ by $\nu_N(v)$. 

\begin{restatable}{claim}{fiveteenthclaim}\label{claim:adjacent_vertices_of_degree_6}
If $[x,y]$ is a digon and both $x$ and $y$ have degree $6$, then either
\begin{enumerate}[label=(\roman*)]
    \item the two neighbours of $y$  distinct from $x$ have degree at least $8$, or
    \item the two neighbours of $x$ distinct from $y$ have 
       degree at least $8$ and $\nu_N(x)\geq 4$.
\end{enumerate}
\end{restatable}

\begin{proofclaim}
Let $[x,y]$ be a digon in $D$ with $d(x)=d(y)=6$.
By Claim~\ref{claim:degree6_3_or_6_neighbours} $n(x)=n(y)=3$.
Let $u$ and $v$ be the two neighbours of $x$ different from $y$.
By Claim~\ref{claim:degree7_7neighbours}, $u$ and $v$ have degree $6$ or at least $8$.

If $u$ and $v$ are linked by a digon, then
by Claim~\ref{claim:triangle_two_vertices_degree_6_dont_exist}, $u$ and $v$ do not have degree $6$,
so they have degree $8$. Moreover $\nu(u) \geq 2$ and $\nu(v)\geq 2$. Thus $\nu_N(x) = \nu(u)+\nu(v) \geq 4$ and (ii) holds.
Henceforth, we may assume that $u$ and $v$ are not linked by a digon.

Let $D'$ the digraph obtained by removing $x$ and $y$ and identifying
$u$ and $v$ into a single vertex $u \star v$.
We claim that $D'$ is not $3$-dicolourable.
To see that, suppose for contradiction that there exists a $3$-dicolouring $\phi$
of $D'$. Then set $\phi(u)=\phi(v)=\phi(u \star v)$, choose
$\phi(y)$ in $\{1,2,3\} \setminus \phi(N(y)\setminus \{x\})$, and finally choose $\phi(x) $ in $\{1,2,3\} \setminus \{\phi(u \star v),\phi(y)\}$. One can easily see that $\phi$
is now a $3$-dicolouring of $D$, a contradiction.
This proves that $D'$ is not $3$-dicolourable
and so it contains a $4$-dicritical digraph $\Tilde{D}$, which must contain $u\star v$ because every subdigraph of $D$ is $3$-dicolourable.
Let $R$ be the subdigraph of $D$ induced by $(V(\Tilde{D}) \setminus \{u \star v\})  \cup  \{u,v,x\}$. We have
\begin{itemize}
    \item $n(R) = n(\tilde{D}) +2$,
    \item $m(R) \geq m(\tilde{D}) +4$ and
    \item $T(R ) \geq  T(\Tilde{D}-u\star v) +1 \geq T(\Tilde{D})$ because $[x,u]$ is a digon, and by Lemma~\ref{lemma:computation_T_minus_one_vertex}.
\end{itemize}
If $\Tilde{D}$ is not $4$-Ore, then 
$\rho(\Tilde{D}) \leq 1$ by minimality of $D$, and so
\[
{
\everymath={\displaystyle}
\renewcommand{\arraystretch}{2}
\begin{array}{r l l}
    \rho(R) &\leq \rho(\Tilde{D}) + 2\left(\frac{10}{3} + \epsilon\right) - 4  & \\
    &\leq 1 + \frac{8}{3} + 2\epsilon & \\
    &< \rho(D) + 3 + 3 \epsilon - 3 \delta &\text{ because $\epsilon - 3\delta \geq -\frac{1}{3}$.}
\end{array}
}
\]
Similarly, if $\Tilde{D}$ is $4$-Ore but not $\bid{K_4}$, then
\[
{
\everymath={\displaystyle}
\renewcommand{\arraystretch}{2}
\begin{array}{r l l}
    \rho(R) &\leq \rho(\Tilde{D}) + 2\left(\frac{10}{3} + \epsilon\right) - 4  & \\
    &\leq \left(\frac{4}{3}+7\epsilon-4\delta\right) + \frac{8}{3} + 2\epsilon &\text{ by Lemma~\ref{lemma:pot_4ore},} \\
    &= 1 + 3 + 9 \epsilon - 4 \delta &\\
    &< \rho(D) +  3 + 9 \epsilon - 4 \delta \\
    &\leq \rho(D) + 3 + 3 \epsilon - 3 \delta &\text{ because $\delta \geq 6\epsilon$.}
\end{array}
}
\]
In both cases (that is when $\Tilde{D}$ is not $\bid{K_4}$), by Claim~\ref{claim:pot_subgraph_final}, 
$D-R$ is a single vertex of degree $6$, namely $y$.
Then every neighbour $w$ of $y$ different from $x$ has degree at least
$6$ in $\Tilde{D}$ (because $\Tilde{D}$ is $3$-dicritical) and so
has degree at least $8$ in $D$ and (i) holds.

Assume now that $\Tilde{D}$ is a copy of $\bid{K_4}$. Let us denote by $a,b,c$
the vertices of $\Tilde{D}$ different from $u\star v$.
Suppose for a contradiction that $u$ has degree $6$. Then $u$ has exactly three neighbours by Claim~\ref{claim:degree6_3_or_6_neighbours}. If $|N(u)\cap \{a,b,c\}| = 2$, then $D\langle \{u,a,b,c\}\rangle$ is a copy of $\bid{K_4}$ minus a digon, contradicting Claim~\ref{claim:no_collapsible}. If $|N(u) \cap \{a,b,c\}| \leq 1$, then $v$ must be adjacent to at least two vertices of $\{a,b,c\}$ with a digon, and so $D\langle \{v,a,b,c\}\rangle$ contains a copy of $\bid{K_4}$ minus a digon, contradicting Claim~\ref{claim:no_collapsible}.
Hence $u$ has degree at least $8$, and by symmetry so does $v$.
Moreover $D\langle\{a,b,c\}\rangle$ is a bidirected triangle, and so by Claim~\ref{claim:triangle_two_vertices_degree_6_dont_exist},
at least two of these vertices have degree at least $8$
(remember that vertices of degree $7$ are in no digon by Claim~\ref{claim:degree7_7neighbours}).
Hence at least four arcs between $\{u,v\}$ and $\{a,b,c\}$ are
incident to two vertices of degree at least $8$.
In other word, $\nu_N(x) =\nu(u) + \nu(v) \geq 4$, so (ii) holds.
\end{proofclaim}

\begin{restatable}{claim}{sixteenthclaim}\label{claim:connected_components_D6}
Let $C$ be a connected component of $D_6$. Then $C$ is one of the following (see Figure~\ref{fig:structure_V6}):
    \begin{enumerate}[label=(\roman*)]
        \item a single vertex, or
        \item a bidirected path on two vertices, or
        \item a bidirected path on three vertices, whose extremities have
            neighbourhood valency at least~$4$, or
        \item a star on four vertices, whose non-central vertices
            have neighbourhood valency at least $4$.
    \end{enumerate}
\end{restatable}
\begin{figure}
  \begin{minipage}{\linewidth}
    \begin{center}	
      \begin{tikzpicture}[thick,scale=1, every node/.style={transform shape}]
        \tikzset{vertex/.style = {circle,fill=black,minimum size=5pt, inner sep=0pt}}
        \tikzset{edge/.style = {->,> = latex}}
        
        \node[vertex] (x) at (0,0) {};
      	\node[vertex,shift={(4,0)}] (u) at (0,0.5) {};
      	\node[vertex,shift={(4,0)}] (v) at (0,-0.5) {};
      	\draw[edge,bend left=20] (u) to (v);
      	\draw[edge,bend left=20] (v) to (u);
      	
      	\node[vertex,shift={(8,0)}] (a) at (0,0) {};
      	\node[vertex,shift={(8,0)}, label=above:$\nu_N \geq 4$] (b) at (0,1) {};
      	\node[vertex,shift={(8,0)}, label=below:$\nu_N \geq 4$] (c) at (0,-1) {};
      	\draw[edge,bend left=20] (a) to (b);
      	\draw[edge,bend left=20] (b) to (a);
      	\draw[edge,bend left=20] (c) to (a);
      	\draw[edge,bend left=20] (a) to (c);
      	
      	\node[vertex,shift={(12,-0.4)}] (d) at (0,0) {};
      	\node[vertex,shift={(12,-0.4)}, label=above:$\nu_N \geq 4$] (e) at (90:1.3) {};
      	\node[vertex,shift={(12,-0.4)}, label=below:$\nu_N \geq 4$] (f) at (-150:1.3) {};
      	\node[vertex,shift={(12,-0.4)}, label=below:$\nu_N \geq 4$] (g) at (-30:1.3) {};
      	\draw[edge,bend left=20] (d) to (e);
      	\draw[edge,bend left=20] (e) to (d);
      	\draw[edge,bend left=20] (d) to (f);
      	\draw[edge,bend left=20] (f) to (d);
      	\draw[edge,bend left=20] (d) to (g);
      	\draw[edge,bend left=20] (g) to (d);
      	
      \end{tikzpicture}
      \label{fig:structure_V6}
      \caption{The possible connected components of $D_6$.}
    \end{center}    
  \end{minipage}
\end{figure}
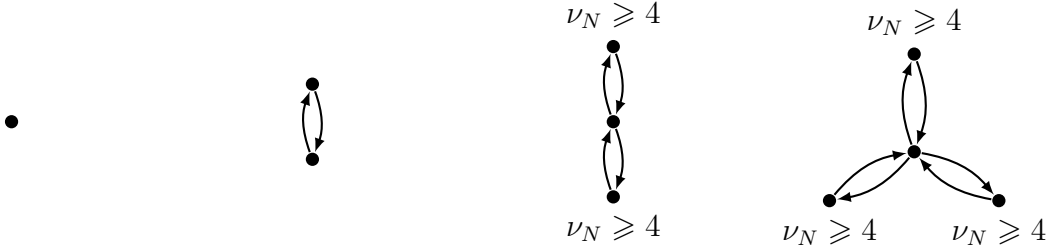
\begin{proofclaim}
First observe that $C$ does not contain a bidirected path $\llbracket x,y,z,w \rrbracket$ on four vertices,
because otherwise, by Claim~\ref{claim:adjacent_vertices_of_degree_6}
applied on $[y,z]$, either $y$ or $z$ has two neighbours of degree at least $8$,
a contradiction.
Observe also that $C$ contains no bidirected triangle by Claim~\ref{claim:triangle_two_vertices_degree_6_dont_exist}.

Moreover, if $\llbracket x,y,z \rrbracket$ is a bidirected path in $C$ on three vertices, then
by Claim~\ref{claim:adjacent_vertices_of_degree_6} applied both on $[y,z]$ and
$[z,y]$, $x$ and $z$ have both neighbourhood valency at least $4$.
The statement of the claim follows.
\end{proofclaim}

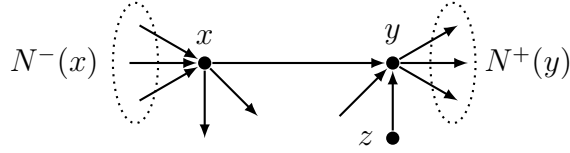
\begin{figure}
  \begin{minipage}{\linewidth}
    \begin{center}	
      \begin{tikzpicture}[thick, every node/.style={transform shape}]
        \tikzset{vertex/.style = {circle,fill=black,minimum size=5pt,
                                        inner sep=0pt}}
        \tikzset{edge/.style = {->,> = latex}}
        \node[vertex, label=above:$x$] (x) at  (0,0) {};
        \node[vertex, label=above:$y$] (y) at  (2.5,0) {};
        \node[vertex, label=left:$z$] (z) at  (2.5,-1) {};
        \draw[edge] (x) to (y);
        \draw[edge] (z) to (y);
        \draw[edge] (x) to (0,-1);
        \draw[edge] (x) to (-45:1);
        \draw[edge] (-1,0) to (x);
        \node[] (Nmx) at (-2,0) {$N^-(x)$};
        \draw[dotted] (-0.92,0) ellipse (0.3 and 0.8);
        \draw[edge] (-150:1) to (x);
        \draw[edge] (150:1) to (x);
        \draw[edge] (y)++(-135:1) to (y);
        \draw[edge] (y) to (2.5+0.866,0.5);
        \draw[edge] (y) to (3.5,0);
        \draw[edge] (y) to (2.5+0.866,-0.5);
        \node[] (Npy) at (4.3,0) {$N^+(y)$};
        \draw[dotted] (3.30,0) ellipse (0.3 and 0.8);
      \end{tikzpicture}
      \caption{An example of an out-chelou arc $xy$.}
      \label{fig:outchelou}
    \end{center}    
  \end{minipage}
\end{figure}

An arc $xy$ is said to be {\bf out-chelou} if
\begin{enumerate}[label=(\roman*)]
    \item $yx \not\in A(D)$,
    \item $d^+(x)=3$,
    \item $d^-(y)=3$, and
    \item there exists $z \in N^-(y) \setminus N^+(y)$ distinct from $x$.
\end{enumerate}
Symmetrically, we say that an arc $xy$ is {\bf in-chelou} if $yx$ is out-chelou in the digraph obtained from $D$ by reversing every arc. See Figure~\ref{fig:outchelou} for an example of an out-chelou arc.

\begin{restatable}{claim}{seventeenthclaim}\label{claim:no_chelou}
    There is no out-chelou arc and no in-chelou arc in $D$.
\end{restatable}
\begin{proofclaim}
By directional duality, it suffices to prove that $D$ has no out-chelou arcs.

Let $xy$ be an out-chelou arc with $z \in N^-(y) \setminus (N^+(y) \cup\{x\})$.
Consider $D' = D - \{x,y\}  \cup \{zz' \mid z' \in N^+(y) \setminus N^-(y)\}$.
We claim that $D'$ is not $3$-dicolourable.
To see that, suppose for contradiction that there is a $3$-dicolouring $\phi$ 
of $D'$.
As $d^+(x)=3$, we can choose $\phi(x)$ in $\{1,2,3\}\setminus \phi(N^+(x)\setminus\{y\})$ to obtain a $3$-dicolouring of $D-y$.
If $|\phi(N^-(y))| < 3$, then choosing $\phi(y)$ in $\{1,2,3\}\setminus \phi(N^-(y))$ gives a $3$-dicolouring of $D$, a contradiction.
Hence $|\phi(N^-(y)| = 3$. Set $\phi(x) = \phi(z)$. Suppose there is a monochromatic directed cycle $C$ in $D$. It must contain $y$ and thus $z$, its unique in-neighbour with its colour. Let $z'$ be the out-neighbour of
$y$ in $C$. It must be in $N^+(y)\setminus N^-(y)$, so $zz'$ is an arc in $D'$. Thus $C - y \cup zz'$ is a monochromatic directed cycle in $D'$,  a contradiction.
Therefore $\phi$ is a $3$-dicolouring of $D$,  a contradiction.
Hence $D'$ is not $3$-dicolourable.

Consequently, $D'$ contains a $4$-dicritical digraph $\Tilde{D}$, which is smaller than $D$ and contains $z$, for otherwise $\Tilde{D}$ would be a subdigraph of $D$.
Consider $U =D\langle V(\Tilde{D}) \cup \{y\}\rangle$.
We have 
\begin{itemize}
    \item $n(U) = n(\Tilde{D}) + 1$,
    \item $m(U) \geq m(\Tilde{D}) + 1$ and
    \item $T(U) \geq T(\Tilde{D}-z)\geq T(\Tilde{D}) -1$ by Lemma~\ref{lemma:computation_T_minus_one_vertex}. 
\end{itemize}  
First if $\Tilde{D}$ is not 4-Ore, then by minimality of $D$ we have $\rho(\Tilde{D}) \leq 1$, so
\[
\rho(U) \leq \rho(\Tilde{D}) + \frac{10}{3}+\epsilon-1+\delta
\leq \frac{10}{3} + \epsilon + \delta \leq \frac{11}{3} - \epsilon - \delta
\]
because $2\epsilon + 2\delta \leq \frac{1}{3}$.

Next if $\Tilde{D}$ is 4-Ore, but not isomorphic to $\bid{K_4}$, then $\rho(\Tilde{D}) \leq \frac{4}{3}+7\epsilon-4\delta$ by Lemma~\ref{lemma:pot_4ore}, and
\[
\rho(U) \leq \rho(\Tilde{D}) + \frac{10}{3}+\epsilon-1+\delta
\leq \frac{11}{3} + 8\epsilon - 3\delta \leq \frac{11}{3} - \epsilon - \delta
\]
because $9\epsilon - 2\delta \leq 0$.

Finally if $\Tilde{D}$ is isomorphic to $\bid{K_4}$, then we have $T(U) \geq T(\Tilde{D}-z)\geq T(\Tilde{D})$ and $\rho(\Tilde{D}) = \frac{4}{3}+4\epsilon-2\delta$.
So the same computation yields
\[
\rho(U) \leq \rho(\Tilde{D}) + \frac{10}{3}+\epsilon-1
\leq \frac{11}{3} + 5\epsilon - 2\delta \leq \frac{11}{3} - \epsilon - \delta
\]
because $6\epsilon - \delta \leq 0$.

In all cases, we have
$\rho(U) \leq \frac{11}{3} - \epsilon - \delta$.
This contradicts Claim~\ref{claim:non_collapsible_pot_bound} because $U$ is not collapsible by Claim~\ref{claim:no_collapsible}.
\end{proofclaim}

\medskip

We now use the discharging method.
For every vertex $v$, let $\sigma(v)=\frac{\delta}{|C|}$ if $v$ has degree $6$ and is in a component $C$ of $D_6$ of size at least $2$, and  $\sigma(v)=0$ otherwise. 
Clearly $T(D)$ is at least the number of connected components of size at least $2$ of $D_6$ so $\sum_{v\in V(D)} \sigma(v) \leq \delta T(D)$.
We define the {\bf initial charge} of $v$ to be $w(v) = \frac{10}{3}+\epsilon-\frac{d(v)}{2}-\sigma(v)$.
We have
\[
\rho(D) \leq \sum_{v\in V(D)} w(v).
\]

We now redistribute this total charge according to the following rules:
\begin{itemize}
    \item[(R1)] A vertex of degree $6$ incident to no digon sends
        $\frac{1}{12}-\frac{\epsilon}{8}$ to each of its neighbours.
    \item[(R2)] A vertex of degree $6$ incident to digons sends
        $\frac{2}{d(v)-\nu(v)}(-\frac{10}{3}+\frac{d(v)}{2}-\epsilon)$ to each neighbour $v$ of degree at least $8$ (so $\frac{1}{d(v)-\nu(v)}(-\frac{10}{3}+\frac{d(v)}{2}-\epsilon)$ via each arc of the digon).
    \item[(R3)] A vertex of degree $7$ with $d^-(v)=3$ (resp. $d^+(v)=3$)
        sends $\frac{1}{12}-\frac{\epsilon}{8}$ to each of its in-neighbours
        (resp. out-neighbours).
\end{itemize}

For every vertex $v$, let $w^*(v)$ be the final charge of $v$.

\begin{restatable}{claim}{eighteenthclaim}
    If $v$ has degree at least $8$, then $w^*(v) \leq 0$.
\end{restatable}
\begin{proofclaim}
Let $v$ be a vertex of degree at least $8$. If $v$ is not adjacent to a vertex of degree at most 7, 
then $w^*(v) = w(v) = \frac{10}{3} + \epsilon - \frac{d(v)}{2} \leq 0$ (because $\epsilon \leq \frac{2}{3}$). 
Otherwise, $d(v) - \nu (v) \geq 1$ and
\begin{align*}
    \frac{1}{d(v)-\nu(v)}\left(-\frac{10}{3}+\frac{d(v)}{2}-\epsilon\right) &\geq \frac{1}{d(v)}\left(-\frac{10}{3}+\frac{d(v)}{2}-\epsilon\right) \\
    &\geq \frac{1}{12} - \frac{\epsilon}{8}.
\end{align*}
Thus $v$ receives at most $\frac{1}{d(v)-\nu(v)}(-\frac{10}{3}+\frac{d(v)}{2}-\epsilon)$
per arc incident with a vertex of degree $6$ or $7$. Since there are $d(v)-\nu(v)$ such arcs,
$w^*(v) \leq w(v) -\frac{10}{3}-\epsilon+\frac{d(v)}{2} = 0$.
\end{proofclaim}

\begin{restatable}{claim}{nineteenthclaim}
    If $v$ has degree $7$, then $w^*(v) \leq 0$.
\end{restatable}
\begin{proofclaim}
By Claim~\ref{claim:degree7_7neighbours}, $v$ has seven neighbours.
Without loss of generality, let us suppose that $d^-(v)=3$ and $d^+(v)=4$.
By Claim~\ref{claim:no_chelou}, the in-neighbours of $v$ can not have out-degree $3$. In particular, they do not have degree $6$, and if they have degree $7$, they do not send anything to $v$ by Rule (R3). 
Hence $v$ receives at most four times
the charge $\frac{1}{12}-\frac{\epsilon}{8}$ by (R1) or (R3), and it sends three times this charge by (R3). 
Hence 
\begin{align*}
    w^*(v) &\leq w(v) +\frac{1}{12}-\frac{\epsilon}{8}\\
            &=-\frac{1}{12}+\frac{7}{8}\epsilon
\end{align*} 
and the result comes because $\epsilon \leq \frac{2}{21}$.
\end{proofclaim}

\begin{restatable}{claim}{twentythclaim}
    If $v$ is a vertex of degree $6$ incident to no digon, then $w^*(v) \leq 0$.
\end{restatable}
\begin{proofclaim}
    The vertex $v$ sends $\frac{1}{12}-\frac{\epsilon}{8}$ to each of its neighbours, and it receives no charge
    as all its in-neighbours (resp. out-neighbours) have out-degree (resp. in-degree) at least $4$, by Claim~\ref{claim:no_chelou}. 
    As a consequence, 
    \[
        w^*(v)  = w(v) -6\left(\frac{1}{12}-\frac{\epsilon}{8}\right) = -\frac{1}{6} +\frac{7\epsilon}{4}
    \]
    and the result comes because $\epsilon \leq \frac{2}{21}$.
\end{proofclaim}

\begin{restatable}{claim}{twentyfirstclaim}\label{claim:computation_charge_D6}
Let $v$ be a vertex in $D_6$ having at least two neighbours of degree
at least $8$. Then $w^*(v) \leq 0$. Moreover, if $v$ is not an isolated vertex in $D_6$ and $\nu_N(v) \geq 4$, then
$w^*(v) \leq -\frac{1}{9}+\frac{5}{3}\epsilon-\frac{\delta}{4}$.
\end{restatable}

\begin{proofclaim}
Observe that $v$ receives no charge and sends the following charge to each of its neighbour $u$ with degree at least $8$:
\begin{align*}
    \frac{2}{d(u)-\nu(u)}\left(-\frac{10}{3}-\epsilon+\frac{d(u)}{2}\right)
        &\geq \frac{2}{d(u)}\left(-\frac{10}{3}-\epsilon+\frac{d(u)}{2}\right)\\
        &= 1 - \frac{2}{d(u)}\left(\frac{10}{3}+\epsilon\right)\\
        &\geq \frac{2}{8}\left(-\frac{10}{3}-\epsilon+4\right)\\
        &= \frac{1}{6}-\frac{\epsilon}{4}.
\end{align*}

\medskip

Assume first that $v$ is isolated in $D_6$. By Claim~\ref{claim:degree7_7neighbours}, its three neighbours do not have degree 7, and so have degree at least 8. Thus $v$ sends three times at least
$\frac{1}{6}-\frac{\epsilon}{4}$, and so 
\[
    w^*(v) \leq w(v) - 3\left( \frac{1}{6} - \frac{\epsilon}{4} \right) = -\frac{1}{6} + \frac{7}{4}\epsilon
\] 
and the result comes because $\epsilon \leq \frac{2}{21}$.

\medskip

Assume now that $v$ is in a connected component $C$ of $D_6$ of size at least $2$.
By Claim~\ref{claim:connected_components_D6}, $\sigma(v) \geq
\frac{\delta}{4}$, so $w(v) \leq \frac{1}{3} +\epsilon -\frac{\delta}{4}$. Moreover, it sends twice at least $\frac{1}{6}-\frac{\epsilon}{4}$.
Hence 
\[
    w^*(v) \leq   \left(\frac{1}{3}+\epsilon-\frac{\delta}{4}\right)-2\left(\frac{1}{6}
-\frac{\epsilon}{4}\right) = \frac{3}{2}\epsilon-\frac{\delta}{4}
\]
and the result comes because $\delta \geq 6\epsilon$. This shows the first part of the statement.

\medskip

We will now prove the second part of the statement. Assume that $v$ is not an isolated vertex in $D_6$ and $\nu_N(v) \geq 4$. Let $u_1$ and $u_2$ be the two neighbours of $v$ with degree at least 8. For every $i\in \{1,2\}$ we have
\[
    \frac{2}{d(u_i)-\nu(u_i)}\left(-\frac{10}{3}-\epsilon+\frac{d(u_i)}{2}\right) = 1 - \frac{1}{d(u_i)-\nu(u_i)}\left(\frac{20}{3}+2\epsilon-\nu(u_i)\right)
\]

\medskip

{\noindent \bf Case 1:}
$\nu(u_i) \geq 7$ for some $i \in \{1,2\}$. Without loss of generality suppose $i=1$.
Then we have
\[
    1 - \frac{1}{d(u_1)-\nu(u_1)}\left(\frac{20}{3}+2\epsilon-\nu(u_1)\right) \geq 1
\]
because $\nu(u_1) \geq 7 \geq \frac{20}{3}+2\epsilon$ as $\epsilon \leq \frac{1}{6}$. Then the total charge sent by $v$ is at least $1$, and thus 
\[
    w^*(v) \leq w(v) -1 \leq \left( \frac{1}{3} + \epsilon - \frac{\delta}{4}\right) - 1 = -\frac{2}{3} + \epsilon - \frac{\delta}{4}
\]
Thus, we have $w^*(v) \leq -\frac{1}{9}+\frac{5}{3}\epsilon-\frac{\delta}{4}$ because $\epsilon, \delta \geq 0$.

\medskip

{\noindent \bf Case 2:} $\nu(u_1),\nu(u_2) \leq 6$.
Let $f: [0,6] \to \mathbb{R}$ be the function defined by 
\[ f(x) = \frac{2}{8-x}\left(-\frac{10}{3}-\epsilon+\frac{8}{2} \right) = 1-\frac{1}{8-x}\left( \frac{20}{3}-2\epsilon-x \right)\]
for every $x \in [0,6]$.
Observe that $f$ is non-decreasing and convex on $[0,6]$ because $-\frac{10}{3}-\epsilon+\frac{8}{2} \geq 0$.
For $i=1,2$, we have 
\[
\frac{2}{d(u_i)-\nu(u_i)}\left(-\frac{10}{3}-\epsilon+\frac{d(u_i)}{2}\right) \geq f(\nu(u_i))
\]
because the function $d \mapsto 1 - \frac{1}{d-\nu(u_i)}\left(\frac{20}{3}+2\epsilon-\nu(u_i)\right)$ is non-decreasing on $[8,+\infty[$ as $\nu(u_i) \leq 6 \leq \frac{20}{3}+2\epsilon$.
Hence the charge sent by $v$ to $u_i$ is at least $f(\nu(u_i))$.
By hypothesis we have $\nu_N(v) = \nu(u_1)+\nu(u_2) \geq 4$.
It follows that the total charge sent by $v$ is at least
\[
{
\everymath={\displaystyle}
\renewcommand{\arraystretch}{2}
\begin{array}{r l l}
    f(\nu(u_1))+f(\nu(u_2)) &\geq 2f\left(\frac{\nu(u_1)+\nu(u_2)}{2}\right) &\text{by convexity of $f$} \\
    &\geq 2f(2)  &\text{because $f$ is non-decreasing} \\
    &= \frac{4}{9}-\frac{2}{3}\epsilon. \\
\end{array}
}
\]
Hence 
\[
    w^*(v) \leq w(v) - \left( \frac{4}{9} - \frac{2}{3}\epsilon \right) \leq \left( \frac{1}{3}+\epsilon-\frac{\delta}{4} \right) -\frac{4}{9}+\frac{2}{3}\epsilon 
= -\frac{1}{9}+\frac{5}{3}\epsilon-\frac{\delta}{4}.
\] 
showing the second part of the statement.
\end{proofclaim}

\begin{restatable}{claim}{twentysecondclaim}
    If $C$ is a connected component of $D_6$, then $\sum_{v \in V(C)}w^*(v) \leq 0$.
\end{restatable}
\begin{proofclaim}
If $C$ has a unique vertex $v$, then, by Claim~\ref{claim:computation_charge_D6},
we have $w^*(v) \leq 0$ as wanted.

If $C$ has two vertices $x$ and $y$, then, again by Claim~\ref{claim:computation_charge_D6},
$w^*(x),w^*(y) \leq 0$, and so $w^*(x)+w^*(y) \leq 0$.

If $C$ is a bidirected path $[x,y,z]$, then,
by Claim~\ref{claim:connected_components_D6}, $x$ and $z$ have both 
neighbourhood valency at least $4$ and so
by Claim~\ref{claim:computation_charge_D6}
$w^*(x),w^*(z) \leq -\frac{1}{9}-\frac{\epsilon}{6}$.
Moreover, $y$ sends at least $\frac{2}{8}(-\frac{10}{3}+4-\epsilon)=
\frac{1}{6}-\frac{\epsilon}{4}$ to its neighbour out of $C$.
Hence 
\[
w^*(y) \leq w(y)- \left(\frac{1}{6}-\frac{\epsilon}{4}\right) \leq  \frac{1}{3}+\epsilon - \frac{\delta}{3}-\frac{1}{6}+\frac{\epsilon}{4} 
    = \frac{1}{6}+\frac{5}{4}\epsilon - \frac{\delta}{3}.
\]
Altogether, we get that
\[
    w^*(x)+w^*(y)+w^*(z) 
    \leq \frac{1}{6}+\frac{5}{4}\epsilon - \frac{\delta}{3} + 2\left(- \frac{1}{9}-\frac{\epsilon}{6} \right)
    =-\frac{1}{18}+\frac{11}{12}\epsilon - \frac{\delta}{3} \leq 0
\]
because $\delta \geq 6\epsilon$.

Finally, if $C$ is a bidirected star with centre $x$ and three other vertices $y,z,w$, then $w^*(x)\leq w(x) = \frac{1}{3}+\epsilon - \frac{\delta}{4}$. Moreover, each of $y,z,w$ has 
neighbourhood valency at least $4$ by Claim~\ref{claim:connected_components_D6}
and so has final charge at most $-\frac{1}{9}+\frac{5}{3}\epsilon-\frac{\delta}{4}$
by Claim~\ref{claim:computation_charge_D6}.
Hence
\[
    w^*(x)+w^*(y)+w^*(z)+w^*(w) 
\leq \frac{1}{3}+\epsilon-\frac{\delta}{4} + 3\left( -\frac{1}{9}+\frac{5}{3}\epsilon-\frac{\delta}{4} \right) \leq 6\epsilon - \delta \leq 0
\]
because $\delta \geq 6\epsilon$.
\end{proofclaim}

As a consequence of these last claims, we have 
$\rho(D) \leq \sum_{v \in V(D)} w(v)=\sum_{v \in V(D)} w^*(v) \leq 0 \leq 1$, a contradiction. 
This proves Theorem~\ref{thm:main_thm_potential}. \qed

\section{An upper bound on \texorpdfstring{$o_k(n)$}{ok(n)}}
\label{section:upperbounds}

In this section, we show that, for every fixed $k$, there are infinitely many values of $n$ such that $o_k(n)\leq (2k-\frac{7}{2})n$. The proof is strongly based on the proof of~\cite[Theorem 4.4]{aboulker3dicritical}, which shows $o_k(n)\leq (2k-3)n$ for every $k,n$ (with $n$ large enough).
For $k=4$, the construction implies in particular that there is a $4$-dicritical oriented graph with $76$ vertices and $330$ arcs, and there are infinitely many $4$-dicritical oriented graphs with $m/n \leq 9/2$.

\begin{proposition}
    Let $k\geq 3$ be an integer. For infinitely many values of $n\in \mathbb{N}$, there exists a $k$-dicritical oriented graph $\vec{G}_k$ on $n$ vertices with at most $(2k-\frac{7}{2})n$ arcs.
\end{proposition}
\begin{proof}
    Let us fix $n_0 \in \mathbb{N}$. We will show, by induction on $k$, that there exists a $k$-dicritical oriented graph $\vec{G}_k$ on $n$ vertices with at most $(2k-\frac{7}{2})n$ arcs, such that $n\geq n_0$.

    \medskip
    
    When $k=3$, the result is known (\cite[Corollary 4.3]{aboulker3dicritical}). We briefly describe the construction for completeness. Start from any orientation of an odd cycle on $2n_0+1$ vertices. Then for each arc $xy$ in this orientation, add a directed triangle $\Vec{C_3}$ and every arc from $y$ to $V(\Vec{C_3})$ and every arc from $V(\Vec{C_3})$ to $x$ (see Figure~\ref{fig:3dicritic}). This gadget forces $x$ and $y$ to have different colours in every $2$-dicolouring. Since we started from an orientation of an odd cycle, the result is a 3-dicritical oriented graph on $4(2n_0+1)$ vertices and $10(2n_0+1)$ arcs.

    \begin{figure}[hbtp]
        \begin{center}	
          \begin{tikzpicture}[thick,scale=1, every node/.style={transform shape}]
            \tikzset{edge/.style = {->,> = latex}}
            \tikzset{vertex/.style = {circle,fill=black,minimum size=5pt, inner sep=0pt}}
            \tikzset{bigvertex/.style = {shape=circle,draw, minimum size=2em}}
            \node[vertex] (3) at  (180:1) {};
            \node[vertex] (2) at  (108:1) {};
            \node[vertex] (1) at  (36:1) {};
            \node[vertex] (0) at  (-36:1) {};
            \node[vertex] (4) at  (-108:1) {};
            \foreach \i in {0,...,2}{
              \pgfmathtruncatemacro{\j}{\i + 1}
              \draw[edge] (\i) to (\j){};
          	}
            \draw[edge] (4) to (3);
            \draw[edge] (4) to (0);
            \node[bigvertex] (u) at  (72:1.8) {$\vec{C_3}$};
            \node[bigvertex] (v) at  (144:1.8) {$\vec{C_3}$};
            \node[bigvertex] (w) at  (-144:1.8) {$\vec{C_3}$};
            \node[bigvertex] (x) at  (-72:1.8) {$\vec{C_3}$};
            \node[bigvertex] (y) at  (0:1.8) {$\vec{C_3}$};
            \draw[edge] (2) to (u);		
    	\draw[edge] (u) to (1);
            \draw[edge] (3) to (v);		
    	\draw[edge] (v) to (2);
            \draw[edge] (w) to (4);		
    	\draw[edge] (3) to (w);
            \draw[edge] (0) to (x);		
    	\draw[edge] (x) to (4);
            \draw[edge] (1) to (y);		
    	\draw[edge] (y) to (0);
          \end{tikzpicture}
          \caption{A 3-dicritical oriented graph with $\frac{5}{2}n$ arcs.}
          \label{fig:3dicritic}
        \end{center}
    \end{figure}
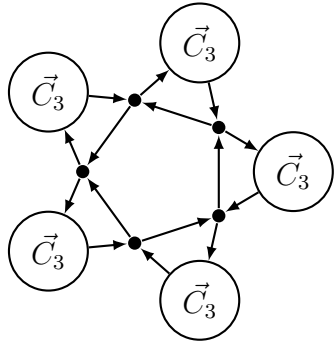

    \medskip 
    
    Let us fix $k\geq 4$ and assume that there exists such a $(k-1)$-dicritical oriented graph $\Vec{G}_{k-1}$ on $n_{k-1} \geq n_0$ vertices with $m_{k-1} \leq (2(k-1)-\frac{7}{2})n_{k-1}$ arcs. We start from any tournament $T$ on $k$ vertices. Then we add, for each arc $xy$ of $T$, a copy $\vec{G}_{k-1}^{xy}$ of $\vec{G}_{k-1}$, all arcs from $y$ to $\vec{G}_{k-1}^{xy}$ and all arcs from $\vec{G}_{k-1}^{xy}$ to $x$. Figure~\ref{fig:4dicritic} illustrates a possible construction of $\vec{G}_4$, where $T$ is the transitive tournament on $4$ vertices.

    \begin{figure}[hbtp]
        \begin{center}	
          \begin{tikzpicture}[thick,scale=1, every node/.style={transform shape}]
            \tikzset{edge/.style = {->,> = latex}}
            \tikzset{vertex/.style = {circle,fill=black,minimum size=5pt, inner sep=0pt}}
            \tikzset{bigvertex/.style = {shape=circle,draw, minimum size=2em}}
            \node[vertex] (0) at  (45:0.8) {};
            \node[vertex] (1) at  (135:0.8) {};
            \node[vertex] (2) at  (-135:0.8) {};
            \node[vertex] (3) at  (-45:0.8) {};
            \foreach \i in {1,...,3}{
              \draw[edge] (0) to (\i){};
            }
            \foreach \i in {2,...,3}{
              \draw[edge] (1) to (\i){};
      	}
            \node[bigvertex] (a) at  (0:2) {$\vec{G}_3$};
            \node[bigvertex] (b) at  (180:2) {$\vec{G}_3$};
            \node[bigvertex] (c) at  (120:2) {$\vec{G}_3$};
            \node[bigvertex] (d) at  (-120:2) {$\vec{G}_3$};
            \node[bigvertex] (e) at  (-60:2) {$\vec{G}_3$};
            \node[bigvertex] (f) at  (60:2) {$\vec{G}_3$};
            \draw[edge] (2) to (3);
            \draw[edge] (2) to (b);
            \draw[edge] (b) to (1);
            \draw[edge] (3) to (a);
            \draw[edge] (a) to (0);
            \draw[edge] (1) to (c);
            \draw[edge] (c) to (0);
            \draw[edge] (3) to (d);
            \draw[edge] (d) to (2);
            \draw[edge] (2) to (e);
            \draw[edge] (e) to (0);
            \draw[edge] (3) to (f);
            \draw[edge] (f) to (1);
          \end{tikzpicture}
          \caption{A 4-dicritical oriented graph with at most $\frac{9}{2}n$ arcs.}
          \label{fig:4dicritic}
        \end{center}
    \end{figure}
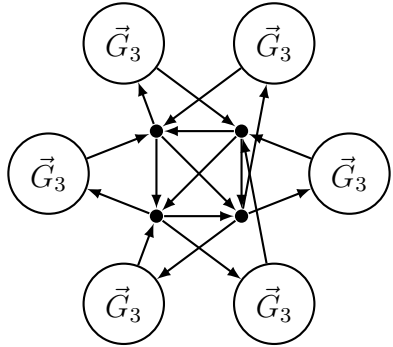

    Let $\vec{G}_k$ be the resulting oriented graph. By construction, $n_k = |V(\Vec{G}_k)|$ and $m_k = |A(\Vec{G}_k)|$ satisfy:
    \begin{align*}
        n_k &= k + \binom{k}{2} n_{k-1}\\
        m_k &= \binom{k}{2} + \binom{k}{2} \times 2 \times n_{k-1} + \binom{k}{2} \times m_{k-1} \\
        &\leq \binom{k}{2} + \binom{k}{2} \left(2 + 2(k-1)-\frac{7}{2}\right)n_{k-1}\\
        &= \binom{k}{2} + \binom{k}{2} \left(2k-\frac{7}{2}\right)n_{k-1}\\
        &= \binom{k}{2} + \left(2k-\frac{7}{2}\right)(n_{k}-k)\\
        &\leq \left(2k-\frac{7}{2}\right)n_{k}
    \end{align*}
    where in the last inequality we used $k\left( 2k-\frac{7}{2}\right) \geq \binom{k}{2}$, which holds when $k\geq 2$.   
    We will now prove that $\Vec{G}_k$ is indeed $k$-dicritical. 
    
    \medskip
    
    We first prove that $\dic(\Vec{G}_k) = k$. Assume that there exists a $(k-1)$-dicolouring $\alpha$ of $\Vec{G}_k$. Then there exist $x,y\in V(T)$ such that $\alpha(x) = \alpha(y)$. Since $\dic(\Vec{G}_{k-1}) = k-1$, there exists $z\in V(\Vec{G}_{k-1}^{xy})$ such that $\alpha(z) = \alpha(x)$. But then $(x,y,z,x)$ is a monochromatic directed triangle in $\alpha$: a contradiction.

    \medskip

    Let us now prove that $\dic(\Vec{G}_{k} \setminus \{uv\}) \leq k-1$ for every arc
    $uv \in A(\Vec{G}_{k})$. This implies immediately that $\dic(\Vec{G}_{k} =k$ and shows the result.
    
    Consider first an arc $uv$ in $A(T)$. We colour each copy $\Vec{G}_{k-1}^{xy}$ of $\Vec{G}_{k-1}$ with a $(k-1)$-dicolouring of $\Vec{G}_{k-1}$. We then choose a distinct colour for every vertex in $T$, except $u$ and $v$ which receive the same colour. This results in a $(k-1)$-dicolouring of $\Vec{G}_{k} \setminus \{uv\}$.

    Consider now an arc $uv$ of $\Vec{G}_{k-1}^{xy}$ for some $xy\in A(T)$. Because $\Vec{G}_{k-1}$ is $(k-1)$-dicritical,  there exists a $(k-2)$-dicolouring $\xi$ of $\Vec{G}_{k-1}^{xy} \setminus \{uv\}$. Hence we colour $\Vec{G}_{k-1}^{xy}\setminus \{uv\}$ with $\xi$, every other copy of $\Vec{G}_{k-1}$ a $(k-1)$-dicolouring of $\Vec{G}_{k-1}$, and we choose a distinct colour for every vertex in $T$, except $x$ and $y$ which both receive colour $k-1$. This results in a $(k-1)$-dicolouring of $\Vec{G}_{k} \setminus \{uv\}$.

    Consider finally an arc $uv$ arc from $u \in V(T)$ to $v\in V(\Vec{G}_{k-1}^{uy})$ (the case of $u \in V(\Vec{G}_{k-1}^{xv})$ and $v\in V(T)$ being symmetric). Because $\Vec{G}_{k-1}$ is dicritical,  there exists a $(k-1)$-dicolouring $\gamma$ of $\Vec{G}_{k-1}^{uy}$ in which $v$ is the only vertex coloured $k-1$. Hence, we colour $\Vec{G}_{k-1}^{uy}$ with $\gamma$, every other copy of $\Vec{G}_{k-1}$ with a $(k-1)$-dicolouring of $\Vec{G}_{k-1}$, and we choose a distinct colour for every vertex in $T$, except $u$ and $y$ which both receive colour $k-1$. This results in a $(k-1)$-dicolouring of $\Vec{G}_{k} \setminus \{uv\}$.
\end{proof}

\bibliographystyle{plain}
\bibliography{biblio}

\end{document}